\documentclass[12pt,a4paper]{amsart} 

\usepackage{stmaryrd}
\usepackage{DLdef1}
\makeatletter
\theoremstyle{thm}
\@ifundefined{Theorem}{}{}
\theoremstyle{definition}
\@ifundefined{Remark}{}{}
\@ifundefined{Remarks}{}{}
\@ifundefined{Example}{\newtheorem*{Example}{\secno Example}}{}
\makeatother
\providecommand{\fG}{\mathfrak{G}}
\providecommand{\what}{\widehat}

\usepackage{hyperref}

\usepackage{bbm}

\newcommand {\fbr}   {{\mathfrak{br}}}
\newcommand {\fer}  {{\mathfrak{er}}}

\newcommand {\ffr}   {{\mathfrak{fr}}}
\newcommand {\fy}   {{\mathfrak{y}}}
\newcommand {\fdy}   {{\mathfrak{dy}}}
\newcommand {\fby}   {{\mathfrak{by}}}
\newcommand {\fsby}   {{\mathfrak{sby}}}
\newcommand {\fme}   {{\mathfrak{me}}}
\newcommand {\fmy}   {{\mathfrak{my}}}
\newcommand {\fsmy}   {{\mathfrak{smy}}}

\newcommand{\un}{\underline{N}}
\newcommand{\del}{\partial}

\newcommand {\Lra} {\Longrightarrow}

\begin{document}

\title[$G(2)$-type structures]{Structures of $G(2)$ type and
non-integrable distributions in characteristic $p$}

\author{Pavel Grozman${}^a$, Dimitry Leites${}^b$}

\address{${}^a$deceased\\
${}^b$Department of Mathematics, Stockholm University, Albanov\"agen 28,
SE-114 19 Stockholm, Sweden; dimleites@gmail.com}

\keywords {Cartan prolongation, non-holonomic manifold,
$G(2)$-structure; Melikyan algebra, Brown algebra, Ermolaev algebra, Frank
algebra, Skryabin algebra}

\subjclass[2020]{Primary 17B50, 17B20; Secondary 70F25}

\begin{abstract} Lately we observe: (1) an upsurge of interest
(in particular, triggered by a~paper by Atiyah and Witten) to
manifolds with $G(2)$-type structure; (2) classifications are
obtained of simple (finite-dimensional and graded vectorial) Lie
superalgebras over fields of complex and real numbers and of simple
finite-dimensional Lie algebras over algebraically closed fields of
characteristic $p$ greater than 3; (3) importance of non-integrable
distributions in observations (1) and (2).

We add to interrelation of (1)--(3) an explicit description of
several exceptional simple Lie algebras for $p=2, 3$ (Brown, Ermolaev,
Frank, and Skryabin algebras, and analogs of Melikyan algebras) as
subalgebras of Lie algebras of vector fields preserving
non-integrable distributions analogous to (or identical with) those
preserved by $\mathfrak{g}(2)$, $\mathfrak{o}(7)$,
$\mathfrak{sp}(4)$, $\mathfrak{sp}(10)$ and the Brown algebra
$\mathfrak{br}(3)$. The description is performed in terms of
Cartan-Tanaka-Shchepochkina prolongs and is similar to descriptions
of simple Lie superalgebras of vector fields with polynomial
coefficients. Our results illustrate usefulness of Shchepochkina's
algorithm and SuperLie package; at least two families of simple Lie
algebras found in the process are new.
\end{abstract}

\maketitle

\begin{quote}\rightline{In memory of Felix Aleksandrovich Berezin}\end{quote}

\thispagestyle{empty}

\setcounter{tocdepth}{3}
\tableofcontents

\
\section{Introduction}

\subsection{What is new in this version} This is an improved version
of the paper published in 2005, see \cite{GL}. In this version we dare to formulate
\textbf{Conjecture: for $p=3$, all simple finite-dimensional Lie
algebras are the examples described in this paper, plus those obtained by
means of the Kostrikin-Shafarevich approach, and the deforms\footnote{\textit{Results} of
deformations, the term coined by M.~Gerstenhaber by analogy with
transform and construct, the results of the corresponding actions.}
of the results of both procedures}.

This version corrects Theorem \ref{Poincare} (although a~mere
illustration given ``to widen the stock of our experience" that does
not affect main results, it should better be correctly formulated).
The material is rearranged and edited; exposition is clearer.  

We
denote by $\Par\un$ the \textit{number of parameters the shearing vector
$\un$ depends on} and resolve the mystery of the concealed parameter
of the shearing vector of the big Skryabin algebra $\fby(\un)$ left
open in the published version of this paper. The allegedly 1-parameter deform of 
$\fh_{\exp}(2;\underline{(1,n)})$ is exorcised. 

References are updated.

Recently we learned about two results published in not easily
accessible sources (\cite{KuJa, KD}). An important improvement of
the formulation of the method due Kostrikin and Shafarevich, due to
Dzhumadildaev and Kostrikin (see \cite{KD}), claims that --- ``DK-claim" --- \textbf{all
simple finite-dimensional Lie algebras over $\Kee$ are obtained
either by the KSh-method or as \textit{deforms}
of certain ``standard" examples}. This improvement definitely works
for $p\geq 5$. In particular, for $p=5$, it removes the Melikyan
algebras from the list of ``standard" examples, see \cite{KD}. In
\cite{MeZu} this claim is double-checked for the particular
case\footnote{This case is very interesting for various reasons,
see \cite{BLLS}.} $\un=\One$, where
\begin{equation}
\label{N_s}
\One:=(1, \dots, 1).
\end{equation}

Kuznetsov and Jakovlev identified the simple Ermolaev algebra, which
we denote\footnote{In this paper, the superscript ${}^{(1)}$ denotes
the first derived algebra or commutant; the prime ${}'$ is used for other
purposes.} $\fer^{(1)}(\One)$, see Subsection~ \ref{fer}, as a~deform
of the contact Lie algebra $\fk^{(1)}(3;\One)$, see the claim in
\cite{KuJa}, corroborated in \cite{GKL}. 
(To prove the DK-claim is still an
\textbf{open problem}. For other open problems, see Subsection~
\ref{Qs}.) Here we investigate which of the simple Lie algebras we
have studied, or even introduced, qualify for the role of
``standard" examples.

In this version, we pay more attention to ``concealed parameters" of
the shearing vectors, parameters not visible in one interpretation
might become manifest in another one.

\subsection{An overview of the stock of simple Lie algebras}
In what follows the ground field $\Kee$ of characteristic $p>0$ is
assumed to be algebraically closed although we use this
property of $\Kee$ only in eq.~\eqref{o18}. 
\textbf{Our results are obtained, mostly, with the
aid of Grozman's \textit{SuperLie} package} now supported by Andrey Krutov, see \cite{GK}.

\sssec{Over $\Cee$} The exceptional (non--serial) simple finite-dimensional Lie algebras, although of considerable interest lately
on account of various applications (\cite{AW, B, CJ, FG}) are far
less understood than, say, $\fsl(n)$, cf. \cite{A, B}. Who, experts
in representation theory including, can nowadays lucidly explain
what is\footnote{We denote the exceptional Lie algebras in the same
way as the serial ones, like $\fsl(n)$; we thus avoid confusing
$\fg(2)$ with the second component $\fg_2$ of a~$\Zee$-graded Lie
algebra $\fg=\oplus \fg_i$. We follow Bourbaki's convention to
denote Lie groups and algebraic groups by Latin characters starting with a
capital letter using same but small characters in Fraktur font for
their Lie algebras. In characteristic $p>0$, the researchers customarily used
Latin capitals for Lie algebras as well.} $\fg(2)$ or $\ff(4)$ or
$\fe(6)-\fe(8)$?! Definitions in terms of octonions, although
beautiful (\cite{BE}), is one way to understand these algebras.
Descriptions in terms of defining relations (as in \cite{GL1}) are
satisfactory for computers, not humans. Here we offer one more way
to interpret certain simple Lie algebras and their ``relatives", such
as their central extensions and algebras of derivations.

This paper gives some applications of
Shchepochkina's general algorithm \cite{Shch} describing Lie
algebras and Lie superalgebras in terms of non--integrable
distributions they preserve. Berezin, who taught all three of us (PG, DL and IShch),
liked to read classics and advised his students to follow his
example.
We return to Cartan's first, now practically forgotten, description
of Lie algebras, not necessarily simple or exceptional ones, in
terms of non--integrable distributions they preserve; we intend to
apply it to ever wider range and begin with $\fo(7)$, $\fsp(4)$ and
$\fsp(10)$.

\sssec{Over fields $\Kee$} We encounter more of \textit{seemingly}
``strange\rq\rq exa\-mples. In \cite[V.1, Introduction]{S}, Strade listed all simple
finite-dimensional Lie algebras over $\Kee$ for $p>3$ and examples
for $p=3$ known to him at that time. We mainly use notation from \cite[v.1]{S} except
for Skryabin algebras which we endow with adjectives. We find
Skryabin's own notation $Y$ with its inherent implicit question mark
most appropriate for ``mysterious" algebras; we only converted Latin
$Y$ to Fraktur font $\fy$: we follow Bourbaki,
leaving capital Latin letters for Lie groups and algebraic groups.

Skryabin discovered several new simple Lie algebras \cite{Sk3,Sk1,Sk2};
Kuznetsov considered various cases of classification of simple finite-dimensional Lie
algebras for $p=3$, see \cite{Ku}.

For deformations, and on classification for $p=2$, see \cite{BGL1,
BGLLS, BLLS}.

Melikyan algebras, still described as something somewhat mysterious
and usually only for $p=5$, are, as we will see, no more mysterious
than $\fg(2)$ for which Shchepochkina recalls Cartan's lucid
description, see \cite{Shch}. We observe that, for $p=2$ and $3$,
Melikyan algebras are the conventional special vectorial (divergence
free) Lie algebras. A more appropriate version of Melikyan algebras
for $p=2$ is the prolong of one of Shen algebras --- a~``correct
version of $\fg(2)$ for $p=2$ --- considered by Brown \cite{Br}; for
clarifications, see \cite{BGLLS2}.

\subsection{Lie algebras preserving non--integrable distributions are vectorial Lie algebras
realized as generalized Cartan prolongs} A.~Kostrikin and
Shafarevich used flags in description of simple Lie algebras in
characteristic $p>0$. In these descriptions, (non--integrable)
distributions appear twice: as associated with flags, and with Lie
algebras of depth $>1$.

Kostrikin felt the
importance of the Cartan prolong and its generalization (for the definitions, see Subsection~
\ref{CTS}, \ref{CTSS}) to algebras
of depth $>1$ but his voice, even amplified by the authority of an
ICM talk, was not heard, except by Elsting, Ermolaev and Kuznetsov
who buried their results in a~little-known journal \textit{Izvestiya Vuzov}.
Examples Kostrikin and his students unearthed (for example, Melikyan
and Ermolaev algebras, followed by Skryabin algebras), as well as
Kuznetsov's interpretations --- practically identical to ours --- of
several of the known algebras, are all obtained as such
(generalized) Cartan prolongs. Still, no general definition of
generalized prolongs --- a~most vital tool --- was ever formulated
for $p>0$ (except a~tentative one in \cite{FSh}; the paper \cite{Shch}
positively answers questions of \cite{FSh}); a~similarity between
these examples and Shchepochkina's constructions of simple Lie
superalgebras, as well as ``non-standard" regradings, was
mentioned only in \cite{KL} and never before or after until recently
this method helped us to discover new serial and exceptional simple Lie (super)algebras,
see \cite{BGLLS2} and refs therein.
This fact and
the lack of lucid algorithm for constructing generalized prolongs
were the reasons why the examples we consider
had to be elucidated and interpreted.


Remarkably, an interpretation we have in mind --- \textbf{the
description of (the exceptional simple) Lie algebras as preserving a
non--integrable distribution and, perhaps, something else} --- WAS
repeatedly published; first, by Cartan, cf. \cite{C}. At least, for
$p=0$. But this aspect of \cite{Y}, as well as of \cite{C}, passed
unnoticed. In \cite{Y}, Yamaguchi lucidly described Tanaka's
construction of generalized prolongs and considered, among other
interesting things, two of the three possible $\Zee$-gradings of
$\fg(2)$ related with ``selected\rq\rq (see eq. (\ref{selec}))
simple roots and interpreted $\fg(2)$ as preserving the
non--integrable distributions associated with these $\Zee$-gradings.

The initial Cartan's interpretation of $\fg(2)$ used one of these
distributions without indication why this particular distribution
was selected. Later Cartan considered another distribution which
characterizes the Hilbert-Cartan equation $f'=(g'')^2$, see \cite{Y} and its recent superization \cite{KST}. Larsson
\cite{La2} considered the remaining, third grading of $\fg(2)$ and
several (selected randomly, it seems) gradings of depth $\leq 2$ of
$\ff(4)$ and $\fe(6)$--$\fe(8)$. These and similar results for other
algebras looked as \textit{ad hoc} examples.

Shchepochkina's algorithm (\cite{Shch}) describes Lie algebras and
superalgebras $\fg$ of vector fields as \textbf{ge\-ne\-ralized}
Cartan prolongs and \textbf{partial} prolongs (\cite{LSh}) in terms
of non--integrable distributions $\fg$ preserves and is applicable to
fields of prime characteristic. Such an interpretation (except
partial prolongs) was known to the classics (Lie, \'E.~Cartan) but
an \textit{explicit} description of the Lie (super) algebras in
terms of non--integrable distributions they preserve was only obtained
(as far as we know) for \textit{some} of the ``selected\rq\rq
gradings of some algebras. We believe it is time for a~thorough
study of all possible distributions related with simple Lie
(super)algebras, and start with \cite{GL3}.

With Shchepochkina's algorithm we immediately see that various
examples previously somewhat mysterious, e.g., Frank algebras, are
just \textit{partial prolongs} corresponding to analogs of the
projective embedding $\fsl(n+1)\subset\fvect(n)$. Likewise, $\fg(2)$
is a~partial prolong if $p=5$ or 2, whereas the Melikyan algebras
are complete prolongs; $\fg(2)$ is the complete prolong only if
$p\neq 5, 2$.

\subsection{Main results. Computer-aided study} For $p=3$, the two
Brown algebras, and their deformations, were until now given by
means of Cartan matrices $A$ with only implicit defining relations
(\ref{myst}), see below. As we pass to the Skryabin algebras, the
mist thickens so much that only non--positive part of one
$\Zee$-graded simple Lie algebra (we denote it $\fby(7;\un)^{(1)}$)
was known, cf. \cite{Sk3} \cite[V.1, Introduction]{S}. Here we determine $\Par\un$ for the Lie
algebras considered; in particular, we discover ``concealed"
parameters for $\fby(7;\un)$ and confirm that $\Par\un$ for the Ermolaev algebra
$\fer(3;\un)$ was determined correctly in earlier papers.

Although certain deforms of the simple derived of Hamiltonian Lie
algebras were known (in particular, Skryabin described all
\textbf{filtered} deforms, see \cite{Sk2}), there was no
classification of deforms, and hence $\fbr(2; 0, 1, c;
\underline{(11n)})$ is a~new simple Lie algebra, same as $\fsby(7;\un
)^{(1)}$.

In what follows we describe the main ideas that lead to these discoveries. Actually,
all these ideas were known, apart from certain technical details developed in \cite{Shch},
but the importance of these ideas was never appreciated as much as they deserve.

We suggest to consider the following simple Lie algebras as
``standard". First of all, $\fbr(3)$. If the notion of a~``standard"
algebra is considered as a~``building brick" for obtaining the other
simple Lie algebras by means of CTS-prolongation (see Subsection~
\ref{CTSS}) of the non--positive part of this ``brick" and subsequent
deforms thereof, then we may stop here. If we take a~less broad
point of view, and allow to deform the ``standard" algebras but not
apply CTS-prolongs to their parts of degree $\leq 0$ or partial
prolongs to their parts of degree $\leq 1$, we have the following
``standard" examples, all but $\fbr(3)$, $\fer^{(1)}(3;\un)$ and
$\ffr^{(1)}(3;\un)$ discovered by Skryabin:
\begin{equation}\label{standAlgs}
\scriptsize
\begin{array}{ll}
\begin{array}{l}\fbr(3);\\
\fby^{(1)}(7;\un), \ \fsby^{(1)}(7;\un), \ \ \Par\un=4 \\
\fdy^{(1)}(10;\un), \ \ \Par\un=3\\\end{array}&\begin{array}{l}
\fmy(6; \un), \ \fsmy^{(1)}(6; \un), \ \ \Par\un=3;\\
\fer^{(1)}(3;\un), \ \ \Par\un=2;\\
\ffr^{(1)}(3;\un), \ \ \Par\un=1.
\end{array}\end{array}
\end{equation}
We found out the correct number of parameters $\Par\un$ the shearing
vector $\un$ depends on for $\fby^{(1)}(7;\un)$, $\fsby^{(1)}(7;\un
)$ and $\fer^{(1)}(3;\un)$.

Since $\dim \fsmy^{(1)}(\One)=\dim \fe^{(1)}(6)/\fc$, it is
tempting to conjecture that one of these algebras is a~deform of the
other one, which is ``standard". In reality, Chebochko proved long time 
ago that $\fe^{(1)}(6)/\fc$ is rigid, see \cite{Ch1}; S.~Bouarroudj et al.
recently proved that no cocycle of non--zero 
weight can represent a~class of $H^2(\fg;\fg)$ for
$\fg=\fsmy^{(1)}(\One)$. Since one algebra is symmetric and the
other one lopsided with respect to its maximal torus, it follows
that no cocycle of weight 0 can deform $\fsmy^{(1)}(\One)$ into
$\fe^{(1)}(6)/\fc$.

Even the algebras that were considered to be known revealed
\textbf{concealed parameters}: in several cases above, not all
parameters on which the shearing vector depends were known. We
observed the following fact:
\begin{equation}\label{factUn}
\begin{minipage}[l]{14cm}
\textbf{If the number of coordinates of $\un$ is greater than
$\Par\un$, then $\Par\un$ is equal to the number of indeterminates
of degree $>1$ in a~certain $\Zee$-grading}.
\end{minipage}
\end{equation}
In examples known to us for $p\leq 5$, \textbf{all} the indeterminates of
degree $>1$ contribute to $\Par\un$, \textbf{except for} $\fdy^{(1)}(10;\un)$ and
$\fme(5;\un)$: for them, $\Par\un$ is equal to the number of
indeterminates of the highest degree in a~certain $\Zee$-grading (the one we consider in this paper), the heights of other indeterminates are constrained (can not exceed~1).

\section{Background: $p=0$}

\ssec{Cartan prolongs}\label{CTS} Let $\fg_0$ be a~Lie algebra,
$\fg_{-1}$ a~$\fg_0$-module. Let us define the $\Zee$-graded Lie
algebra $(\fg_{-1}, \fg_{0})_{*}=\mathop{\oplus}\limits_{i\geq
-1}\fg_i$ called the \textit{complete Cartan prolong} (the result of
the Cartan \textit{prolongation}) of the pair $(\fg_{-1}, \fg_{0})$.
Geometrically the Cartan prolong is the maximal Lie algebra of
symmetries of the $G$-structure (here: $\fg_0=\text{Lie}(G)$) on
$\fg_{-1}$. The components $\fg_i$ for $i>0$ are defined
recursively.

First, recall that, for any (finite-dimensional)
vector space $V$, we have
\begin{equation}\label{3.}
\Hom(V, \Hom(V,\ldots, \Hom(V,V)\ldots)) \simeq L^{i}(V, V, \ldots,
V; V),
\end{equation}
where $L^{i}$ is the space of $i$-linear maps and we have
$(i+1)$-many $V$'s on both sides. Now, we recursively define, for
any $v_1, \dots, v_{i+1}\in \fg_{-1}$ and any $i > 0$:
\begin{equation}\label{4.}
\renewcommand{\arraystretch}{1.4}
\begin{array}{ll}
\fg_i =& \{X\in \Hom(\fg_{-1}, \fg_{i-1})\mid X(v_1)(v_2, v_3,
...,
v_{i+1}) =X(v_2)(v_1, v_3, ..., v_{i+1})\}.
\end{array}
\end{equation}

Let the $\fg_0$-module $\fg_{-1}$ be faithful. Then, clearly,
\begin{equation}\label{5.}
(\fg_{-1}, \fg_{0})_{*}:=\oplus\fg_i\subset \fvect (m) = \fder~
\Cee[[x_1, \ldots , x_m ]], \; \text{ where}\; m = \dim~ \fg_{-1}.
\end{equation}
Moreover, setting $\deg x_i=1$ for all $i$, we see that
\begin{equation}\label{6.}
\renewcommand{\arraystretch}{1.4}
\begin{array}{l}
\fg_i = \{X\in \fvect (m)\mid \deg X=i,\; [X, \del]\in
\fg_{i-1}\;\text{ for any }\del\in \fg_{-1}\}.
\end{array}
\end{equation}

Now it is subject to an easy verification that the Cartan prolong
$(\fg_{-1}, \fg_{0})_{*}$ forms a~subalgebra of $\fvect (n)$. (It is
also easy to see that $(\fg_{-1}, \fg_{0})_{*}$ is a~Lie algebra
even if $\fg_{-1}$ is not a~faithful $\fg_0$-module.)

\ssec{Nonholonomic manifolds. Cartan-Tanaka-Shchepochkina
prolongs}\label{CTSS} Let $M^n$ be an $n$-di\-men\-si\-o\-nal
manifold with a~non--integrable distribution $\cD$. Let
\begin{equation}\label{7.}
\cD= \cD_{-1}\subset \cD_{-2} \subset \cD_{-3} \dots \subset
\cD_{-d}
\end{equation}
be the sequence of strict inclusions, where the fiber of $\cD_{-i}$
at a~point $x\in M$ is
\begin{equation}\label{.8}
\cD_{-i+1}(x)+ [\cD_{-1}, \cD_{-i+1}](x)
\end{equation}
(here $[\cD_{-1}, \cD_{-i-1}]=\Span\left([X, Y]\mid X\in
\Gamma(\cD_{-1}), Y\in \Gamma(\cD_{-i-1})\right)$) and $d$ is the
least number such that
\begin{equation}\label{.9}
\cD_{-d}(x)+[\cD_{-1}, \cD_{-d}](x) = \cD_{-d}(x).
\end{equation}
In case $\cD_{-d} = TM$ the distribution is called
\textit{completely non--holonomic}. The number $d = d(M)$ is called
the \textit{non--holonomicity degree}. A manifold $M$ with a
distribution $\cD$ on it will be referred to as
\textit{non--holonomic} one if $d(M)\neq 1$. Let
\begin{equation}\label{n_i}
n_i(x) = \dim \cD_{-i}(x); \qquad n_0(x)=0; \qquad
n_d(x)=n-n_{d-1}.
\end{equation}
The distribution $\cD$ is said to be \textit{regular} if all the
dimensions $n_i$ are constants on $M$. We will only consider
regular, completely non--holonomic distributions, and, moreover,
satisfying certain transitivity condition (\ref{tr}) introduced
below.

To the tangent bundle over a~non--holonomic manifold $(M, \cD)$ we
assign a~bundle of $\Zee$-graded nilpotent Lie algebras as
follows. Fix a~point $pt\in M$. The usual adic filtration by
powers of the maximal ideal $\fm:=\fm_{pt}$ consisting of
functions that vanish at $pt$ should be modified because distinct
coordinates may have distinct ``degrees". The distribution
$\cD$ induces the following filtration in $\fm$:
\begin{equation}\label{2.1}\scriptsize
\renewcommand{\arraystretch}{1.4}
\begin{array}{ll}
\fm_k=&\{f\in\fm\mid X_1^{a_1}\ldots X_n^{a_n}(f)=0\;\text{ for
any $X_1, \dots,X_{n_1}\in \Gamma(\cD_{-1})$, }\\
& \text{$X_{n_1+1}, \dots, X_{n_2}\in \Gamma(\cD_{-2})$,\dots,
$X_{n_{d-1}+1}, \dots, X_{n}\in \Gamma(\cD_{-d})$}\\
& \text{ such that }\; \mathop{\sum}\limits_{1\leq i\leq d}\quad
i\mathop{\sum}\limits_{n_{i-1}< j\leq n_{i}} a_j\leq k\},
\end{array}
\end{equation}
where $\Gamma(\cD_{-j})$ is the space of germs at $pt$ of sections
of the bundle $\cD_{-j}$. Now, to a~filtration
\begin{equation}\label{.10}
\cD= \cD_{-1}\subset \cD_{-2} \subset \cD_{-3} \dots \subset
\cD_{-d}=TM, \end{equation} we assign the associated graded bundle
\begin{equation}\label{.11}
\gr(TM)=\oplus\gr\cD_{-i},\;\text{ where
$\gr\cD_{-i}=\cD_{-i}/\cD_{-i+1}$}
\end{equation}
and the bracket of sections of $\gr(TM)$ is, by definition, the
one induced by bracketing vector fields, the sections of $TM$. We assume
a ``transitivity condition\rq\rq: The
Lie algebras
\begin{equation}\label{tr}
\gr(TM)|_{pt}
\end{equation}
induced at each point $pt\in M$ are isomorphic.

The grading of the coordinates $(\ref{2.1})$ determines a
non--standard grading of $\fvect(n)$ (recall (\ref{n_i})):
\begin{equation}\label{gr}\scriptsize
\renewcommand{\arraystretch}{1.4}
\begin{array}{l}
\deg x_1=\ldots =\deg x_{n_1}=1,\\
\deg x_{n_1+1}=\ldots =\deg x_{n_2}=2,\\
\dotfill \\
\deg x_{n-n_{d-1}+1}=\ldots =\deg x_{n}=d.
\end{array}
\end{equation}
Denote by $\fv=\mathop{\oplus}\limits_{i\geq -d}\fv_i$ the algebra
$\fvect(n)$ with the grading $(\ref{gr})$. One can show that the
``complete prolong'' of $\fg_-$ to be defined shortly, i.e.,
$(\fg_-)_{*}:=(\fg_-, \widetilde \fg_0)_{*}\subset \fv$, where
$\widetilde \fg_0:=\fder_0\fg_-$, preserves $\cD$.

For non--holonomic manifolds, an analog of the group $G$ from the term
``$G$-structure'', or rather of its Lie algebra,
$\fg=\text{Lie}(G)$, is the pair $(\fg_-, \fg_0)$, where $\fg_0$ is
a subalgebra of the $\Zee$-grading preserving Lie algebra of
derivations of $\fg_-$, i.e., $\fg_0 \subset \fder_0\,\fg_-$. If
$\fg_0$ is not explicitly indicated, we assume that $\fg_0
=\fder_0\,\fg_-$, i.e., is the largest possible.

Given a~pair $(\fg_-, \fg_0)$ as above, define its
\textit{generalized Cartan prolong} or
\textit{Cartan-Tanaka-Shchepochkina prolong} (briefly: CTS-prolong)
to be the maximal subalgebra $(\fg_-,
\fg_0)_{*}=\mathop{\oplus}\limits_{k\geq -d} \fg_k$ of $\fv$ with
given non--positive part $(\fg_-, \fg_0)$. For an explicit
construction of the components, see \cite{Sh14,Y,Shch}.

Since the time the manuscript of this article was written we came to the following short  definition of the prolong. Recall that the graded Lie superalgebra
$\fb=\oplus_{k \geqslant -d}\ \fb_k$
is said to be \textit{transitive}\index{Lie (super)algebra, transitive} if for all~ $k\geq 0$ we have
\begin{equation*}
\label{12}
\{x\in \fb_k \mid [x,\fb_{-}]=0\}=0, \text{~~where $\fb_{-}:=\oplus _{k <0}\,\fb_k$}.
\end{equation*}

Let
$\fg_-=\mathop{\oplus}_{-d\leq i\leq -1}\ \fg_i$ be a~nilpotent
$\Zee$-graded Lie (super)algebra and let $\fg_0$ be a~Lie sub(super)algebra of the Lie
(super)algebra $\mathfrak{der}_0(\fg_-)$ of degree-0 derivations of~ $\fg_-$. 
The maximal transitive $\Zee$-graded Lie (super)algebra $\fg=\oplus_{k \geqslant -d}\ \fg_k$ whose non-positive part is the given $\fg_-\oplus \fg_0$ is called the 
\textit{Cartan prolong}\index{Cartan prolongation} 
--- the result of Cartan prolongation --- of the pair $(\fg_-,\fg_0)$ and is denoted by $(\fg_-,\fg_0)_*$. 

\ssec{Partial prolongs and projective structures} Let $(\fg_-,
\fg_0)_{*}$ be a~depth $d$ Lie algebra; $\fh_1\subset \fg_1$ be a
$\fg_0$-submodule such that $[\fg_{-1}, \fh_1]=\fg_0$. If such
$\fh_1$ exists (usually, $[\fg_{-1}, \fh_1]\subset\fg_0$), define
the $i$th partial prolong of $(\mathop{\oplus}\limits_{i\leq
0}\fg_i, \fh_1)$ for $i\geq 2$ to be
\begin{equation}
\label{partprol} \fh_{i}=\{D\in\fg_{i}\mid [D, \fg_{-1}]\in
\fh_{i-1}\}.
\end{equation}
Set $\fh_i=\fg_i$ for $i\leq 0$ and call
$\fh_{*}=\mathop{\oplus}\limits_{i\geq -d}\fh_i$ the Shchepochkina
\textit{partial prolong} of $(\mathop{\oplus}\limits_{i\leq 0}\fg_i,
\fh_1)$, see \cite{Sh14}. (Of course, the partial prolong can also
be defined if $\fh_0$ is contained in $\fg_0$.)

\begin{Example} The $SL(n+1)$-action on the projective space
$\Pee^n$ gives the embedding 
\[
\fpgl(n+1):=\fgl(n+1)/\Kee 1_{n+1}\simeq\fsl(n+1)\subset\fvect(n); 
\]
here
$\fsl(n+1)$ is a~partial prolong of $\fvect(n)_{i\leq
0}\oplus\fh_1$ for some $\fh_1$.
\end{Example}

\ssec{The two types of Lie algebras} For the list of simple Lie
superalgebras (finite-dimensional and $\Zee$-graded of polynomial
growth) and background on Linear Algebra in Superspaces, see
\cite{LSh, BGLLS2}. All this super knowledge is not a~must to
\textit{understand} this paper, but comparison of the situation over
$\Kee$ and super cases over $\Cee$ is instructive. Observe that
there are only two major types of Lie (super)algebras:

(SY) For \textit{symmetric} algebras, related with a~Cartan
subalgebra is a~root decomposition such that ($\sdim$ is
\textit{superdimension})
\begin{equation}
\label{sy}
 \sdim\fg_{\alpha}=\sdim\fg_{-\alpha}\;\text{ for any root
 }\;\alpha;
\end{equation}

(LOP) For \textit{lopsided} algebras, related with a~Cartan
subalgebra is a~root decomposition such that (\ref{sy}) fails.
(Usually, lopsided algebras can be realized as \textbf{vectorial}
Lie superalgebras.)

\ssec{Integer bases in Lie superalgebras} Let $A=(A_{ij})$ be
an $n\times n$ matrix. A \textit{Lie superalgebra $\fg=\fg(A)$ with
Cartan matrix} $A=(A_{ij})$, is given by its \textit{Chevalley
generators}, i.e., elements $X_i^{\pm}$ of degree $\pm 1$ and $H_i =
[X_{i}^+, X_{i}^-]$ (of degree 0) that satisfy the relations
(hereafter in similar occasions either all superscripts $\pm$ are
$+$ or all are $-$)
\begin{equation}
\label{g(a)} {}[X_{i}^+, X_{j}^-] = \delta_{ij}H_i, \quad [H_{i},
H_{j}] = 0, \quad [H_i, X_j^{\pm}] = \pm A_{ij}X_j^{\pm},
\end{equation}
and additional relations $R_i=0$ whose left sides are implicitly
described, for a~general Cartan matrix, as
\begin{equation}
\label{myst}
\renewcommand{\arraystretch}{1.4}
\begin{array}{l}
\text{``the $R_i$ that generate the maximal ideal $I$ such
that}\\
I\cap\Span(H_i \mid 1\leq i\leq n)=0. "
\end{array}
\end{equation}
For simple (finite-dimensional) Lie algebras over $\Cee$, instead of
implicit description (\ref{myst}) we have the following explicit
description (\textit{Serre relations}): Normalize $A$ so that
$A_{ii}=2$ for all $i$; then the off-diagonal elements of $A$ are
non--positive integers and
\begin{equation}
\label{serre} (\ad_{X_i^{\pm}})^{1-A_{ij}} (X_j^{\pm}) = 0 .
\end{equation}
A way to normalize $A$ may affect reduction modulo $p$: Letting some
diagonal elements of the \textbf{integer} matrix $A$ be equal to 1
we make the Cartan matrices of $\fo(2n+1)$ and Lie superalgebra
$\fosp(1|2n)$ (for definition of Cartan matrices of Lie
superalgebras, see \cite{GL1}) indistinguishable (this accounts for
their ``remarkable likeness'' \cite{RS,Ser,DLZh}):
\begin{equation}
\label{2ways}\scriptsize
\begin{pmatrix}
 \ddots&\ddots&\ddots&\vdots\\
 \ldots & 2 & -1 & 0 \\
 \ldots & -1 & 2 &-2 \\
 \ldots & 0 &-2 & 2
\end{pmatrix} \; \text{ or } \;
\begin{pmatrix}
\ddots&\ddots&\ddots&\vdots\\
\ldots & 2 & -1 & 0 \\
\ldots & -1 & 2 &-2 \\
\ldots & 0 &-1 & 1
\end{pmatrix}
\end{equation}
\textbf{For Lie superalgebras of the form $\fg=\fg(A)$}, there exist
bases with respect to which all structure constants are integer. Up
to the above indicated two ways (\ref{2ways}) to normalize $A$,
there is only one such (\textit{Chevalley}) basis, cf. \cite{Er}.

When $p=3$ and 2, it may happen that $A_{ii}=0$ (if $p=0$ and $p>3$,
then $A_{ii}=0$ implies ${\dim\fg(A)=\infty}$). It is natural to study
this case in terms of vectorial Lie algebras.

\textbf{For vectorial Lie superalgebras}, integer bases are
associated with $\Zee$-forms of $\Cee[x]$ --- a~supercommutative
superalgebra in $a$ (ordered for convenience) indeterminates $x =
(x_1,...,x_a)$ of which the first $m$ indeterminates are even, the
remaining $n$ ones are odd ($m+n=a$). For an $m+k$-tuple of
non--negative integers $\underline{r}=(r_1, \ldots , r_a)$, where
$r_i=0$ or 1 for $i>m$, we set
\begin{equation}
\label{1.} u_i^{(r_{i})} := \frac{x_i^{r_{i}}}{r_i!}\quad
\text{and}\quad u^{(\underline{r})} := \prod\limits_{1\leq i\leq a}
u_i^{(r_{i})}.
\end{equation}
Let us formally replace fractions with $r_i!$ in denominators by
inseparable symbols $u_i^{(r_{i})}$ which are well-defined over
fields of prime characteristic. Clearly,
\begin{equation}
\label{divp} u^{(\underline{r})} \cdot u^{(\underline{s})} = \binom
{\underline{r} + \underline{s}} {\underline{r}} u^{(\underline{r} +
\underline{s})}, \quad \text{where}\quad\binom {\underline{r} +
\underline{s}} {\underline{r}}:=\prod\limits_{1\leq i\leq a}\binom
{r_{i} + s_{i}} {r_{i}}.
\end{equation}


\ssec{Traces and divergences on vectorial Lie algebras} The
\emph{traces} --- various linear functionals that vanish on the
first derived subalgebra, a.k.a. the commutant, $\fg^{(1)}:=[\fg,\fg]$ of $\fg$ ---
belong to one type of analogs of the trace on the matrix Lie
algebras. The \emph{divergences} (depending on a~fixed volume
element) belong to another type. For an interpretation of
divergences as ``prolongations" of traces, see \cite{LeP}.

Accordingly, the \emph{special} or \emph{divergence-free} subalgebra
of a~vectorial algebra $\fg$ is denoted by $\fs\fg$, e.g.,
$\fsvect(n)$; the codimension of $\fg^{(1)}$ in $\fg$
is equal to the number of traces on $\fg$.

\section{Background: $p>0$}

\ssec{Functions and vector fields} The elements of $\Zee/n$ are
denoted by $\bar a$, where $a\in \Zee$ represents the class $\bar a$. For an $m$-tuple of positive
integers $\un = (N_1,..., N_m)$, denote
\begin{equation}
\label{u;N}
\renewcommand{\arraystretch}{1.4}
\begin{array}{l}
\cO(m; \un):=\Kee[u; \un ]:=
\Span_{\Kee}(u^{\underline{(r)}}\mid r_i < p^{N_{i}}\;\text{ for any
$i$}).\end{array}
\end{equation}
As is clear from (\ref{divp}), $\Kee[u; \un ]$ is a
subalgebra of $\Kee[u]$. The algebra $\Kee[u]$ and its subalgebras
$\Kee[u; \un ]$ are called the \textit{algebras of divided
powers}; they are analogs of the algebra of functions.

Since any derivation $D$ of a~given algebra is determined by the
values of $D$ on the generators, we see that the Lie algebra
$\fder(\cO(m; \un))$ has more than $m$ functional
parameters (coefficients of the analogs of partial derivatives) if
$N_i\neq 1$ for at least one $i$. Define \textit{distinguished
partial derivatives} by setting
\begin{equation}
\label{disder} \del_i(u_j^{(k)})=\delta_{ij}u_j^{(k-1)}\;\text{ for all
$k<p^{N_j}$}.
\end{equation}
The \textit{general} vectorial Lie algebra is
\begin{equation}
\label{2.} \fvect(m; \un) \;\text{ a.k.a
 } W(m; \un)\; :=\fder_{dist} \cO(m; \un)=\{\sum f_i\del_i\mid f_i\in\cO(m; \un)\text{~~for $i=1,\dots, m$}\}.
\end{equation}

\sssec{Fock spaces}\label{hei} For Lie superalgebras $\fg=\fg(A)$, if $A_{ii}=0$, then $x^{+}_i$
and $x^{-}_i$ generate an analog of the Heisenberg Lie algebra: if
$x^{\pm}_i$ are even, we denote this Lie algebra $\fhei(2; p;
\un)$, where $N\in \Nee$. Its natural representation is
realized in the Fock space of functions $\cO(1; \un)$; this module
is indecomposable for $\un >1$ and irreducible for
$\un =1$.

If the $x^{\pm}_i$, for a~fixed $i$, are odd, they generate
$\fsl(1|1; \Kee)$; all the non--trivial irreducible $\fsl(1|1; \Kee)$-modules
are of dimension $1|1$.

\ssec{On modules over vectorial Lie algebras} For simple complex
vectorial Lie algebras with polynomial or formal coefficients
considered with a~natural $x$-adic topology (as algebras of
continuous derivatives of $\Cee[[x]]$), Rudakov described the
irreducible continuous representations. Up to dualization (passage
to induced modules), all such representations either depend on
$k$-jets of the vector fields for $k>1$ and are coinduced, or $k=1$
and then the irreducible representations are realized in the spaces
of tensor fields and also are coinduced, except for the spaces
$\Omega^i$ of differential $i$-forms which have submodules
$Z^i:=\{\omega \in\Omega^i\mid d\omega=0\}$. For a~review of results
on classifications of invariant differential operators over $\Cee$, see
\cite{GLS}. The spaces $Z^i$ are irreducible since (Poincar\'e's
lemma) on the $m$-dimensional space with coordinates $x$, the spaces
$\Omega^i$ are united by the exterior differential $d$ into the
following exact sequence
\begin{equation}
\label{ext}
0\tto\Cee\tto\Omega^0\stackrel{d}{\tto}\Omega^1\stackrel{d}{\tto}\Omega^2
\stackrel{d}{\tto}\dots\stackrel{d}{\tto}\Omega^m\tto 0.
\end{equation}

At the moment, there is no complete description of irreducible
representations of simple finite-dimensional vectorial Lie
algebras for $p>0$. For partial results, see works of Ya.~Krylyuk, e.g., \cite{Kr}, and \cite{BjL1}. 
The next theorem first appeared, it seems,
in Krylyuk's Ph.D. thesis and was
cited in \cite{Sk1,Sk2}.

\sssbegin{Theorem}[Krylyuk]\label{Poincare}
Over $\Kee$, let $\vvol(u):=du_1\wedge\dots\wedge du_m\in
\Omega^m(m;\un)$ be the volume element. Denote:
\[
Z^i(m;\un):=\{\omega \in\Omega^i(m;\un)\mid d\omega=0\},\qquad
B^i(m;\un):=\{d\omega \mid \omega\in\Omega^{i-1}(m;\un)\}.
\]

Set
\begin{equation}\label{tau}
\tau(\un):=(p^{N_1}-1, \dots, p^{N_m}-1).
\end{equation}
The sequence
\begin{equation}
\label{ext-p}
\renewcommand{\arraystretch}{1.4}
\begin{array}{l}
0\tto\Kee\tto\Omega^0(\underline{m;N})\stackrel{d}{\tto}
\Omega^1(m;\un)\stackrel{d}{\tto}\Omega^2(m;\un)
\stackrel{d}{\tto}\dots \stackrel{d}{\tto}\Omega^m(m;\un)
\stackrel{\int}{\tto}\Kee, \end{array}
\end{equation}
where the integral is an analog of the Berezin integral over superspaces (or vice versa: the Berezin integral is an analog of this integral for positive characteristic)
\begin{equation}
\label{int-p}
\renewcommand{\arraystretch}{1.4}
\int f(u)\vvol(u)=\text{the coefficient of the term $u^{\tau(\un)}\vvol(u)$},
\end{equation}
is not exact: the space
$H^a(m;\un):=Z^a(m;\un)/B^a(m;\un)$ is
spanned by the elements
\begin{equation}
\label{coh} u_{i_1}^{(\tau(\un)_{i_1})}\dots
u_{i_k}^{(\tau(\un)_{i_k})}du_{i_1}\dots du_{i_k},\text{~~where $a=i_1+\dots +i_k$.}
\end{equation}
The spaces $B^i(m; \un)$ are irreducible.
\end{Theorem}

\begin{proof}Induction on $m$. For $m=1$ this is obvious. \end{proof}

For any $\fvect(m; \un |n)$-$\cO(m;
\un |n)$-bimodule $M$ with the $\fvect(m;
\un)$-action $\rho$, we denote by $M_{A\, \Div}$ a~copy of
$M$ with the affine $\fvect(m; \un |n)$-action given by
\begin{equation}
\label{div}
\begin{array}{l}
\rho_{A\,\Div}(D)(\mu)=\rho(D)(\mu)+A\,\Div(D)(\mu)\\
\text{for any $D\in \fvect(m; \un |n),\; \mu\in M$ and
$A\in \Kee$}.\end{array}
\end{equation}

After Strade, we denote the space $\Vol(m; \un)$ of
volume forms by $\cO(m; \un)_{\Div}$; denote the
subspace of forms with integral 0 by
\begin{equation}\label{.38}
\cO'(m; \un)_{\Div}=\Span(u^{(a)}\vvol(u)\mid a_i<
\tau(\un)_i\text{~~for all $i$}).
\end{equation}

\ssec{$\Zee$-gradings} Recall that every $\Zee$-grading of a
given vectorial algebra is determined by setting $\deg u_i=r_i\in
\Zee$; every $\Zee$-grading of a~given Lie superalgebra $\fg(A)$ is
determined by setting ${\deg X^{\pm}_i=\pm r_i\in \Zee}$. For the Lie
algebras of the form $\fg(A)$, we set
\begin{equation}\label{selec}
\deg X^{\pm}_i=\pm\delta_{i, i_j}\;\text{ for any $i_j$ from a
selected set $\{i_1, \dots , i_k\}$}
\end{equation}
and say that we have ``\textbf{selected}" certain $k$ pairs of
Chevalley generators (or respective nodes of the Dynkin graph).
Yamaguchi's theorem quoted below (Subsection~ \ref{YaTh}) shows that, in the study of Cartan
prolongs defined in sec.~2.5, the only gradings that yield prolongs
distinct from the initial Lie algebra are the ones with $1\leq k\leq
2$ pairs of ``\textbf{selected}" Chevalley generators.

\ssec{Yamaguchi's theorem}\label{YaTh} Let 
${(\fs_-)_{*}=(\fs_-, \fg_0)_{*}}$ be the Cartan prolong with the
maximal possible $\fg_0=\fder_0(\fs_-)$, where $\fs=\mathop{\oplus}\limits_{i\geq
-d}\fs_i$ is a~simple finite-dimensional Lie algebra.

\sssbegin{Theorem}[\cite{Y}] Over $\Cee$, equality $(\fs_-)_{*}=\fs$ holds
almost always. The exceptions (cases where
$\fs=\mathop{\oplus}\limits_{i\geq -d}\fs_i$ is a~partial prolong
in $(\fs_-)_{*}=(\fs_-, \fg_0)_{*}$) are

\emph{1)} $\fs$ with the grading of depth $d=1$ $($in which case
$(\fs_-)_{*}=\fvect(\fs_-^*)$$)$;

\emph{2)} $\fs$ with the grading of depth $d=2$ and
$\dim\fs_{-2}=1$, i.e., with the ``contact'' grading, in which case
$(\fs_-)_{*}=\fk(\fs_-^*)$ $($these cases correspond to ``selection"
of the nodes on the Dynkin graph connected with the node for the
maximal root on the extended graph$)$;

\emph{3)} $\fs$ is either $\fsl(n+1)$ or $\fsp(2n)$ with the grading
determined by ``selecting'' the first and the $i$th of simple
coroots, where $1<i<n$ for $\fsl(n+1)$ and $i=n$ for $\fsp(2n)$.
$($Observe that $d=2$ with $\dim \fs_{-2}>1$ for $\fsl(n+1)$ and $d=3$
for $\fsp(2n)$.$)$

Moreover, the equality $(\fs_-, \fs_0)_{*}=\fs$ also holds almost
always. The cases where the equality fails $($the ones where a
projective action is possible$)$ are $\fsl(n+1)$ or $\fsp(2n)$ with
the grading determined by ``selecting'' only one $($the first$)$
simple coroot. \end{Theorem}

Observe that the Yamaguchi theorem (Subsection~ \ref{YaTh}) is derived
from the classification of simple Lie algebras of polynomial growth
over $\Cee$. In the absence of similar classification over $\Kee$ if
characteristic $p=3$ (resp. 2), we \textbf{conjecture} that
\textsl{all ``standard" simple Lie algebras are the generalized
Cartan-Tanaka-Shchepochkina prolongs of non--positive or negative
parts of Lie algebras with indecomposable Cartan matrix, or simple
subquotients of these Lie algebras corresponding to one (resp. at
most two\footnote{If $p=2$, selecting two pairs of Chevalley
generators yields new simple Lie algebras as compared with the stock
of simple Lie algebras obtained by selecting just one pair of
Chevalley generators, see \cite{BGLLS}.}) pairs of ``selected"
Chevalley generators}.

\textbf{In this paper we consider ONE ($k=1$) pair of ``selected"
Chevalley generators}.

For vectorial algebras, filtrations are more natural than gradings;
the very term ``vectorial" means, actually, that the algebra is
endowed with a~particular (\textit{Weisfeiler}) filtration, see
\cite{LSh}.

\section{Examples from the literature}

\centerline{\fbox{\textbf{In what follows in this paper, $p=3$ except for Melikyan algebras}.}}

In this section, we list several Lie algebras more or less as
described in \cite[v.1, Introduction]{S}; in the next section we give their
interpretations in terms of (partial) prolongs: no version of
Yamaguchi's theorem is yet available for $p>0$. For a~general
algorithm that describes the non-holonomic distributions these Lie
algebras preserve, see \cite{Shch}.

\ssec{Melikyan algebras for $p=5$} For any prime $p>2$, on the space
$\fg_{-1}:=\cO(1; 1)/\Kee\,1$ of ``functions (in one indeterminate
$u$) modulo constants", the antisymmetric bilinear form
\begin{equation}\label{int}
(f, g)=\int fg'du \text{~~(recall definition \eqref{int-p})},
\end{equation}
where ${}'$ is the distinguished derivative, is non-degenerate. Hence $\fvect(1; \underline{1})$ is embedded into $\fsp(p-1)$.
So we can consider the prolong
\begin{equation}\label{.13}
\fg_{*,\un}=\mathop{\oplus}\limits_{i\geq
-2}\fg_i:=(\fk(p; \un)_{-}, \fc \fvect(1; \underline{1}))_{*,\un}\subset \fk(p; \un),
\end{equation}
where
$\fc\fg=\fg\oplus \Kee\ z$ is the trivial central extension of
$\fg$ and $\un$ is not fixed in advance but is determined by the inner properties of the prolong (several of its indeterminates are constrained). This construction resembles Shchepochkina's construction of
some of exceptional simple vectorial Lie superalgebras
\cite{Sh14}. If $p>5$ or some of the coordinates of
$\un$ is $>1$, then either $\fg_1=0$ or  $\fg_{*,\un}$ is a~known simple Lie algebra.

In late 1970s, Melikyan observed \cite{Me} that if $p=5$, then the prolong
\[
(\fk(5;\un)_{-}, \fc \fvect(1;
\underline{1}))_{*,\un}\subset \fk(5;\un)
\]
was a~new (at that time) simple Lie algebra, $\fme(\un)$.
Melikyan's only available publication lacked details: he did not
write for which 5-tuples $\un$ it is possible to
generalize the construction to $\fk(5;\un)$ and the
ground for Melikyan's claim that $\Par\un=2$ was unclear; we verified this by computer.
The following expressions of the basis elements of non--positive degree in terms of the indeterminate $u$ of
$\fvect(1; \underline{1})$ rewritten in terms of elements of
$\fk(5;\un)$, i.e., in terms of the indeterminates
$t; p_1, p_2, q_1, q_2$  are convenient for constructing
the elements of $\fme(\un)$ of positive degree:
\begin{equation}\label{meme}
\renewcommand{\arraystretch}{1.4}
\scriptsize{
\begin{tabular}{|c|c|c|c|}
\hline $\fg_i$:&$\fg_0\simeq \fvect(1;1)\oplus \Kee \ z$&$\fg_{-1}$&$\fg_{-2}$\cr \hline
span-&$p_1^{(2)}\leftrightarrow u^{(4)}\frac{d}{du}, \;
-p_1p_2\leftrightarrow u^{(3)}\frac{d}{du},\;
-p_1q_2+2p_2^{(2)}\leftrightarrow u^{(2)}\frac{d}{du},$&
$p_1\leftrightarrow -u^{(4)},\; p_2\leftrightarrow u^{(3)}, $&$1$\cr
ned
by&$-p_1q_1-2p_2q_2\leftrightarrow u\frac{d}{du},\; -p_2q_1+q_2^{(2)}\leftrightarrow
\frac{d}{du};\ \ t\leftrightarrow z$&$ q_2\leftrightarrow u^{(2)},\; q_1\leftrightarrow
u$&\cr \hline
\end{tabular}\\
}
\end{equation}

The table \eqref{meme} corrects the one published in \cite{GL}; the correction is due to I.~Shchepochkina, and here are the details for verification.
First of all, recall the formula for the contact bracket (with a minus in front of the Poisson bracket):
\begin{equation}\label{kbme}
\{f,g\}_{k.b.}=(2-E)(f)\del_t(g)-\del_t(f)(2-E)(g)-\sum(\del_{p_i}(f)\del_{q_i}(g)- \del_{q_i}(f)\del_{p_i}(g)),
\end{equation}
where $E:=\sum y_i\del_{y_i}$ and the $y_i$ are all the indeterminates except~$t$.
Hence, $\{p_i,q_i\}=-1$. The formula
\[
(f,g):=\int fdg =\text{the coefficient of $u^{(4)}$ in the expansion of $fdg$}
\]
defines a~non-degenerate anti-symmetric bilinear form on $\cO(1;1)/\Kee$:
\[
\begin{array}{l}
(u^{(4)},u)=\int u^{(4)}du=1,\ \ \ \ \
(u^{(3)},u^{(2)})=\int u^{(3)}udu=4=-1.
\end{array}
\]
This implies the description of a basis of $\fme_{-1}$:
$$
p_1\leftrightarrow -u^{(4)}, p_2\leftrightarrow u^{(3)}, q_2\leftrightarrow u^{(2)}, q_1\leftrightarrow u.
$$
Therefore,
 $\del_u \leftrightarrow -p_2q_1+q_2^{(2)}$. The action of  $u\del_u$ on $\cO(1;1)$ in terms of $u$ and $p, q$ is as follows:
$$
\scriptsize
u\del_u: \begin{cases}-u^{(4)}\mapsto -u\cdot u^{(3)}=u^{(4)},\; u^{(3)}\mapsto u\cdot u^{(2)}=3u^{(3)}=-2u^{(3)},\; u^{(2)}\mapsto u\cdot u=2u^{(2)},\; u\mapsto u;\\
p_1\mapsto -p_1,\; p_2\mapsto -2p_2,\; q_2\mapsto 2q_2, \; q_1\mapsto q_1.\end{cases}
$$
Hence,
\[
u\del_u \leftrightarrow -p_1q_1-2p_2q_2.
\]

Similarly,
\[
\begin{array}{l}
u^{(2)}\del_u \leftrightarrow -p_1q_2+2p_2^{(2)},\\
u^{(3)}\del_u\leftrightarrow -p_1p_2, \\
u^{(4)}\del_u\leftrightarrow p_1^{(2)}.\\
\end{array}
\]
Clearly, $\fme_{-2}=\Kee 1$, and the center of $\fme_{0}$ is $z\leftrightarrow t$.

In this realization all elements of degree $-1$ and the corresponding indeterminates look as if they are on equal footing and it is absolutely unclear why $\un$ depends only of 2 parameters, the height of the other 3 indeterminates being constrained. However, they are NOT on equal footing! They would have been if they were interchangeable under the action of $\fme_{0}$, i.e., if $\fme_{0}$ had acted on them symmetrically.

\sssec{Description of $\fme_1$} The generating functions independent of $t$ are as follows
$$
p_1^{(3)}, \; p_1^{(2)}p_2, \; p_1^{(2)}q_2-2p_1p_2^{(2)}, \; p_1^{(2)}q_1+2p_1p_2q_2+p_2^{(3)}, \;
p_1q_2^{(2)}-p_1p_2q_1-2p_2^{(2)}q_2.
$$

The fact that these functions belong to $\fme_1$ is subject to a direct verification. The fact that there are no other functions independent of $t$ follows from the dimension considerations since, as is not difficult to see,
\be\label{fG}
\text{the functions independent of $t$ generate the central extension $\fG$ of $\fh(2; r)$ for $r=(1,0)$}.
\ee

Observe that $\fG$ is not isomorphic to the known central extension --- the Poisson Lie algebra $\fpo(2; r)$ with the height of the indeterminate of degree 0 being constrained, see Subsec.~\ref{H^2}. 
In more detail: first, the quotient of the space of functions independent of $t$ modulo $\Kee\ 1$ generates the Hamiltonian Lie algebra~$\fh(2; r)$. However, the generating function 1 in $\fme$ is not the same as generating function 1 in the Poisson Lie algebra. In the regrading $\fpo(2; r)$ of the Poisson Lie algebra on the space of functions in $p$ and $q$, the depth is equal to 1 and all functions depending only on $p$ --- they span $\fpo(2; r)_{-1}$ --- commute. Whereas the depth of the Melikyan algebra $\fme$ is equal to 2, $\fme_{-1}$ is spanned by functions in $p$, except 1, and they do not commute but yield $1\in\fme_{-2}$.

\paragraph{The central extension $\fG$ of $\fh(2;r)$ in \eqref{fG}}\label{H^2}
Since the algebra acts by zero on its cohomology, it is easy to see that the weight 0 (modulo 5) of the cocycle might occur only in products of elements of negative degree, so the central extension can only increase the depth of the algebra $\fh(2; \un)$, so it suffices to calculate the extensions setting $\un=\One$ to simplify computations.
The space $H^2(\fh(2;\One))$ is spanned by the following cocycles whose indexes are their degrees in the grading $(1,-1)$ and $\what{x}$ denotes the element of $\Pi(\fg^*)$ dual to $x\in\fg$:
\begin{align*}
 c_{-5} = {} & \what{p^{(2)}}\wedge \what{p^{(3)}} - \what{p}\wedge\what{p^{(4)}},
 & c_0 = {}& \what{p}\wedge\what{q},
 & c_{5} = {} & \what{q^{(2)}}\wedge \what{q^{(3)}} - \what{q}\wedge\what{q^{(4)}}. 
\end{align*}
 \normalsize
The cocycle $c_0$ corresponds to the Poisson algebra, the other two cocycles correspond to algebras isomorphic
to~$\fG$.

Denote $G:=-p_2q_1+q_2^{(2)}$ the generating function corresponding to $\del_u$. Then, the functions depending on $t$ are as follows:
 \[\scriptsize
 \begin{array}{l}
 F_{p_1}=tp_1+2p_1p_2q_2-p_2^{(3)}, \\
 F_{p_2}=\{F_{p_1},G\}=tp_2+p_1q_2^{(2)}-2p_1p_2q_1, \\
 F_{q_2}=\{F_{p_2},G\}=-tq_2+p_2q_2^{(2)}+p_2^{(2)}q_1+p_1 q_1q_2,\\
 F_{q_1}=\{F_{q_2},G\}=tq_1+2q_2^{(3)}-p_2q_1q_2-2p_1q_1^{(2)}.
 \end{array}
 \]

 The operators corresponding to the generating functions in $\fme_-$ are as follows:
 \begin{equation}\label{K-}
 K_{p_i}:=p_i\del_t-\del_{q_i} \text{~ and~}  K_{q_i}:=q_i\del_t+\del_{p_i} \text{~ for~} i=1,2, \text{~ and~}K_1:=2\del_t.
 \end{equation}

As is explained in \cite{Shch}, to select the generating functions that belong to $\fme_-$ from all generating functions of  $\fk(5)$, we should consider the differential equations expressed in terms of the operators commuting with the operators (\ref{K-}); these are
\begin{equation}\label{hatK-}
 \widetilde K_{p_i}=p_i\del_t+\del_{q_i} \text{~ and~}  \widetilde K_{q_i}=q_i\del_t-\del_{p_i}\text{~ for~} i=1,2, \text{~ and~}
 \widetilde K_1=\del_t.
 \end{equation}

Observe that all generating functions in $\fme_0$ satisfy, in particular, the following equations:
\begin{equation}\label{eqme}
(\widetilde K_{p_1})^2 (f)=0, \quad \widetilde K_{p_1}\widetilde K_{p_2}(f)=0.
\end{equation}

Observe that these are far from all equations that single out $\fme\subset \fk(5)$, but for our purposes they suffice. Observe that the operators $\widetilde K_{p_i}$ act on functions $f$ independent of $t$ just as $\del_{q_1}$. Hence, on such functions, the equations (\ref{eqme}) turn into
\begin{equation}\label{fpq}
\frac{\del^2 f}{\del q_1^2}=0, \quad \frac{\del^2 f}{\del q_1\del q_2}=0 \Lra f=q_1\vf(p_1,p_2)+\psi(p_1,p_2,q_2).
\end{equation}

If $\deg_t f=k>0$, i.e., $f=t^k g(p_1,p_2,q_2,q_1)+h$, where $\deg_t h<k$, then taking $k$ times the contact bracket of $f$ with $1$, we see that $g\in\fme$, and hence is also of the same form as $f$ in (\ref{fpq}).

\sssec{Kuznetsov's description of $\fme(\protect\un)$} Kuznetsov offered another description of $\fme(\un)$, see
\cite{Ku1}. From Yamaguchi's theorem cited above we know that the
Cartan-Tanaka-Shchepoch\-kina prolong of $(\fg(2)_-, \fg(2)_0)$ (in
any $\Zee$-grading of $\fg(2)$) is isomorphic to $\fg(2)$, at least,
over $\Cee$. There are two $\Zee$-gradings of $\fg(2)$ with
\textbf{one} ``selected\rq\rq generator: one of depth 2 and one of
depth 3. Kuznetsov observed that, for $p=5$, the non-positive parts
of $\fg(2)$ in the grading of depth 3 are isomorphic to the
respective non-positive parts of the Melikyan algebras in one of
their $\Zee$-gradings. Let $V[k]$ be the $\fgl(V)$-module which is
$V=V[1]$ as $\fsl(V)$-module, and let the fixed central element $z=\id_V\in
\fgl(V)$ act on $V[k]$ as $k\ \id_V$. Then,
\begin{equation}\label{melg2}
\renewcommand{\arraystretch}{1.4}
\scriptsize{
\begin{tabular}{|c|c|c|c|}
\hline
$\fg_0$&$\fg_{-1}$&$\fg_{-2}$&$\fg_{-3}$\cr \hline
$\fgl(2)\simeq\fgl(V)$&$V=V[-1]$&$E^2(V)$&$V[-3]$\cr
\hline
\end{tabular}
}
\end{equation}
So it is natural to conjecture that, for $p=5$, Melikyan algebras
$\fme(\un)$ are
\textit{complete}\footnote{The Lie algebra $\fme(\un)$ could be smaller
than what it actually is, the complete prolong.}
Cartan-Tanaka-Shchepochkina prolongs $(\fg(2)_-, \fg(2)_0)_{*,\un}$ of total
symmetries preserving a~non-holonomic structure whereas $\fg(2)$ is a
\textit{projective} type subalgebra in $\fme(\un)$. In
this realization, it remains unclear what are the admissible values
of $\un $.

Kuznetsov \cite{Ku1} gave yet another realization of $\fme(\un)$. As
$\Zee/3$-graded Lie algebra, it can be represented as
\begin{equation}\label{defme}
\fme(\un):=\fg_\ev \oplus\fg_{\bar 1}\oplus\fg_{\bar 2}
\simeq \fvect(2; \un)\oplus \widetilde \fvect(2;
\un)_{2\Div}\oplus \cO(2; \un)_{-2\Div}\ ,
\end{equation}
where $\widetilde \fvect(2; \un)$ is a~copy of $\fvect(2; \un)$
endowed, together with each element, with a~tilde to distinguish
from (the elements of) $\fvect(2;\un)$. Let $v$ be a~short for
$\vvol(u)$; observe that we have the following identifications, where $v:=du_1du_2$:
\begin{equation}\label{ident1}
\renewcommand{\arraystretch}{1.4}\scriptsize
\begin{array}{l}
\text{for any $p$, \; $du_iv^{-1}=\sign(ij)\del_j$\; for any
permutation $(ij)$
of $(12)$}\\
\text{for $p=5$, \; $
v^{-4}=v, \quad du_1du_2v^2=v^3=v^{-2}$.}
\end{array}
\end{equation}
The $\fg_\ev$-action on the $\fg_{\bar i}$ is natural; the
multiplication in $\fme(\un)$ is given by the following formulas
(in line 2 we use $du_iv=\sign(ij)\del_jv^{2}$, see (\ref{ident1}),
and set $[fv^{-2}, \widetilde Dv^{2}]:=-[\widetilde Dv^{2},
fv^{-2}]$):
\begin{equation}\label{multme}\scriptsize
\renewcommand{\arraystretch}{1.4}
\begin{array}{ll}
{}[f_1v^{-2}, f_2v^{-2}]&=
2\left(f_1d(f_2)-f_2d(f_1)\right)v^4\\
&=2\left(f_1\del_1(f_2)-f_2\del_1(f_1)\right)du_1v +
cycle(12)\\
&=2\left(f_2\del_2(f_1)-f_1\del_2(f_2)\right)\tilde\del_1v^2 +
cycle(12);\\
{}[fv^{-2}, \widetilde Dv^{2}]&=fD;\\
{}[\sum f_i\widetilde \del_iv^{2}, \sum g_j\widetilde
\del_jv^{2}]&=[\sum \sign(ik)f_idu_kv,
\sum \sign(jl)g_jdu_lv]\\
&=\left(f_1g_2-f_2g_1\right)du_1du_2v^2=\left(f_1g_2-f_2g_1\right)v^{-2}.
\end{array}
\end{equation}
The \textit{standard} $\Zee$-grading of $\fme(\un)$ is given by setting (\cite{S}):
\begin{equation}\label{grme}
\renewcommand{\arraystretch}{1.4}
\begin{array}{l}
\deg u^{(\underline{r})}\del_i =3|\underline{r}|-3,\quad \deg
u^{(\underline{r})}v^{-2} =3|\underline{r}|-2,\quad
\deg u^{(\underline{r})}\tilde\del_iv^{2} =3|\underline{r}|-1.\\
\end{array}
\end{equation}

\sssec{On concealed parameters of $\protect\un$}\label{ConcPar} The realization
\eqref{defme} allows one to easily compute the dimensions of
$\fme(\un)$ and its homogeneous components, shows that $\un $
depends on at least $2$ parameters but does not preclude more. The
upper bound on the number of independent parameters of $\un $ comes
from the classification and \cite{KD}.

\sssec{Melikyan algebras for $p=3$} Shchepochkina's realization
\cite{Shch} of the non-positive part of $\fme(\un)$, identical to
that of $\fg(2)$ in a~$\Zee$-grading (\ref{realme}), only involves
$\pm 1$ as coefficients in $\fg_-$ and $\pm 1, \pm 2$ in $\fg_0$ and
so invites to study the prolongs $(\fg_-, \fg_0)_{*,\un}$ for $p=3$,
and $(\fg_-)_{*,\un}$ for $p=2$; this is being done.

Another approach is to interpret decomposition (\ref{defme}): In
$\fvect(3; \underline{N, 1})$, consider a~non-standard
$\Zee$-grading:
\begin{equation}\label{non-stme}
\deg u_1=\deg u_2=0;\quad \deg u_3=1\ .
\end{equation}
Let $u=(u_1, u_2)$, $v=u_3$; $\del_i=\del_{u_i}$; $\del=\del_v$.
Then, $\fvect(3; \underline{N, 1})$ can be represented as a~direct
sum of the following spaces non-canon-ically identified with a~pair of
3 copies of $\fvect(2; \un)$-modules, where $\langle T \rangle$
denotes the space spanned by the elements of a~set $T$:
\begin{equation}\label{non-stdec}
\renewcommand{\arraystretch}{1.4}
\begin{array}{l}
\fvect(2; \un)\simeq\langle f_i(u)\del_i \rangle; \quad \fvect(2;
\un)\simeq\langle vf_i(u)\del_i \rangle; \quad
\fvect(2; \un)\simeq\langle v^2f_i(u)\del_i \rangle; \\
\cO(2; \un)\simeq\langle f(u)\del\rangle; \quad \cO(2; \un
)\simeq\langle f(u)v\del\rangle; \quad \cO(2; \un)\simeq\langle
f(u)v^2\del\rangle\ .
\end{array}
\end{equation}
If we recall that $2\equiv -1 \mod 3$, we see that the corresponding
decomposition of $\fsvect(3; \underline{N, 1})$ is of the form
(\ref{defme}):
\begin{equation}\label{non-stdecS}\scriptsize
\renewcommand{\arraystretch}{1.4}
\begin{array}{l}
\fvect(2; \un)\simeq\langle \sum f_i(u)\del_i - \left(\sum
\del_i(f_i(u))\right )v\del\rangle; \\
\fvect(2; \un)_{\Div}\simeq\langle vf_i(u)\del_i - \left(\sum
\del_i(f_i(u))\right) v^2\del\rangle; \\
\cO(2; \un)_{2\Div}\simeq\langle f(u)\del\rangle\ .
\end{array}
\end{equation}

Thus, the Melikyan algebras for $p=3$ are (isomorphic to) $\fsvect(3; \underline{N,
1})$. Having observed this we recalled that Shen \cite{Sh} had
noticed that, for $p=2$, $\fg(2)\simeq \fsvect(3; \underline{1, 1,
1})$.

\sssec{Melikyan algebras for $p=2$} It is also natural to consider
the prolongs of non-positive parts of $\fg(2)$ in its various
$\Zee$-gradings for $p=2$, and Brown \cite{Br} did just it: As
$\Zee/3$-graded Lie algebras, let
\begin{equation}\label{defme2}
L(\un):=\fg_\ev \oplus\fg_{\bar 1}\oplus\fg_{\bar 2}
\simeq \fvect(2; \un)\oplus \cO(2;
\un)_{\Div}\oplus \cO(2; \un)\ .
\end{equation}
The $\fg_\ev$-action on the $\fg_{\bar i}$ is natural; the
multiplication in $L(\un)$ is given by the following
formulas:
\begin{equation}\label{multme2}
\renewcommand{\arraystretch}{1.4}
\begin{array}{lcl}
{}[fv, g]&=
&fH_g;\\
{}[f, g]&=&H_f(g)v,
\end{array}
\end{equation}
where
\begin{equation}\label{.14}
H_f=\pderf{f}{u_1}\del_2+\pderf{f}{u_2}\del_1\text{~~ for any
$f\in\cO(2; \un)$}.
\end{equation}
Define a~$\Zee$-grading of $L(\un)$ by setting
\begin{equation}\label{grme2}
\renewcommand{\arraystretch}{1.4}
\begin{array}{l}
\deg u^{(\underline{r})}\del_i =3|\underline{r}|-3,\quad \deg
u^{(\underline{r})}v =3|\underline{r}|-2,\quad
\deg u^{(\underline{r})} =3|\underline{r}|-4.\\
\end{array}
\end{equation}
Now, set
$\fme(\un)=L(\un)/L(\un)_{-4}$. This
algebra is not simple, because $\cO(2; \un)_{\Div}$ has a
submodule of codimension 1; but $\fme^{(1)}(\un)$ is
simple.

As is easy to see, the non-positive parts of $\fg(2)$ and $\fme(\un
)$ are isomorphic; $\fme(\un)$ is the complete CTS-prolong of this
part.

\ssec{Brown algebras} The Cartan matrices of Brown algebras are
as follows (here, of course, $-1$ is the same as 2 modulo 3, but
$-1$ is more adequate in relations of Serre
type):
\begin{equation}\label{br2a}
\scriptsize
\fbr(2; \eps)\text{~with CM~}
\begin{pmatrix}2 & -1 \\
 -2 & 1-\eps\end{pmatrix},\text{~where $\eps\neq 0$, and $\fbr(2):=\fbr(2; 1)$},
\end{equation}
observe that\footnote{All
deforms of $\fo(5)$ are described in \cite{BLW}; in addition to
Brown algebras there is one more simple Lie algebra among the
deforms; in this paper we consider its prolongs in which it serves as
the 0th component.} $\fbr(2; -1)\simeq\fo(5)\simeq\fsp(4)$;
\begin{equation}\label{br3a}
\scriptsize{
1\fbr(3)\text{~with CM~}
\begin{pmatrix}2&-1&0\\-1&2&-1\\
0&-1&0\end{pmatrix}\text{~and ~} 2\fbr(3)\text{~with CM~}
\begin{pmatrix}2&-1&0\\-2&2&-1\\
0&-1&0\end{pmatrix}.
}
\end{equation}
The reflections (similar to those that generate Weyl groups of
simple Lie algebras over $\Cee$, see \cite{CCLL}) interchange the Cartan
matrices \eqref{br3a} of $\fbr(3)$, whereas  for $\fbr(2; \eps)$ these reflections change the value
of the parameter $ \eps\mapsto\eps'=\eps^{-1}$; in \cite{Ko}, it is proved that:
\begin{equation}\label{brso}
\fbr(2;\eps)\simeq \fbr(2;\eps')\text{~~if and only if $\eps\eps'=1$ for $\eps\neq \eps'$}.
\end{equation}
The Brown algebras are sometimes denoted by $\fbr(2;\alpha)$, where
$\alpha=\frac{1}{\eps-1}$, 
which forces one to consider two types of Cartan matrices, see, e.g., \cite{BGL2}; the parameter $\eps$ is more convenient and natural than~ $\alpha$.

Even when relations (\ref{myst}) were implicit it was known
(\cite{S}) that the dimensions of Brown algebras given by the Cartan
matrices \eqref{br2a} and \eqref{br3a} are equal to 10 and 29,
respectively (assuming the usual rules \eqref{g(a)}, \eqref{myst} of
constructing $\fg(A)$ from $A$). In what follows we give explicit
relations between Chevalley generators of the Brown algebras.

Kostrikin \cite{Ko} described a~3-parameter family containing
$\fbr(2)$, see \eqref{br2a}, and acknowledged that Ruda\-kov was the first to observe
that, if $p=3$, then for \textbf{any} irreducible $\fsl(2)$-module
$V$, the Cartan prolong $\fg:=\oplus \fg_{i}$, where $\fg_{-1}=V$
and $\fg_0=\fsl(2)$ or $\fgl(2)$ whose center acts on $\fg$ as the
grading operator, is a~simple Lie algebra. Nobody, it seems,
published so far exact descriptions of this Cartan prolong
(\ref{Tabcg0}), (\ref{g10}) nor were particular cases studied (see
(\ref{br01c})).

\ssec{Cartan prolongs associated with irreducible $\fsl(2)$-modules
for $p=3$} There are not that many non--trivial irreducible $\fsl(2)$-modules;
all such modules are listed in \cite{RS}: for $p=3$, there is just
one module $V$ of dimension 2 --- the tautalogical one; it brings about $\fh(2;\un):=(V, \fsl(V))_{*,\un}$
and $\fk(3;\un):=(\fhei(2), \fgl(V))_{*,\un})$. There is also a~3-parameter family of
3-dimensional modules $\Tee(a, b, c)$ given by the following matrices, where $a\neq
bc$:
\begin{equation}\label{.15}
\scriptsize
\widetilde X^-=\begin{pmatrix}
0&0&c\\
1&0&0\\
0&1&0
\end{pmatrix}\quad \widetilde H=\begin{pmatrix}
a-bc&0&0\\
0&0&0\\
0&0&-a+bc\end{pmatrix} \quad \widetilde X^+=\begin{pmatrix}
0&a&0\\
0&0&a\\
b&0&0\end{pmatrix}.
\end{equation}

Let $\fbr(2; (a, b, c);\un):=(\Tee(a, b, c), \fgl(2))_{*,\un}$.
Having normalized the matrices \eqref{.15} as follows
\begin{equation}\label{.16}\scriptsize
\widetilde X^+=\begin{cases}\begin{pmatrix}
0&1&0\\
0&0&1\\
\displaystyle\frac{b}{a}&0&0\end{pmatrix}&\text{ if $a\neq 0$ and
then $H=\begin{pmatrix}
1-bc&0&0\\
0&0&0\\
0&0&-1+bc\end{pmatrix}$}\cr
\begin{pmatrix}
0&0&0\\
0&0&0\\
1&0&0\end{pmatrix}&\text{ if $b\neq 0$ and then $H=\begin{pmatrix}
-c&0&0\\
0&0&0\\
0&0&c\end{pmatrix}$}\end{cases} 
\end{equation}
we see that $\fbr(2; (a, b, c);\One)$ depends on at most two
parameters; actually, it depends on one parameter, as proved in
\cite{BLW}. Theorem \ref{ProloTabc} states that $\fbr(2; (a, b,
c);\un)$ may have components of degree $>1$ only if $a=0$ and $b=1$.
In realization by vector fields we have
\begin{equation}
\label{Tabcg0}\scriptsize
\renewcommand{\arraystretch}{1.4}
\begin{array}{ll}
\widetilde X^-=c u_1 \del_3 + u_2 \del_1 + u_3 \del_2, &\widetilde
H=(a - b
c) \left(u_1\del_1 -u_3 \del_3\right), \\
\widetilde X^+=a u_1 \del_2 + a~u_2 \del_3 + b u_3 \del_1,
&\widetilde E=u_1 \del_1 + u_2 \del_2+ u_3 \del_3.
\end{array}
\end{equation}
Indeed,
\begin{equation}\label{.17}
{}[\widetilde X^+, \widetilde X^-]=\widetilde H,\qquad [\widetilde
H, \widetilde X^{\pm}]=\pm(a-bc)\widetilde X^{\pm}.
\end{equation}
If $a\neq bc$, the change
\begin{equation}\label{o18}
\widetilde X^{\pm}\mapsto X^{\pm}:=\sqrt{a-bc}\widetilde
X^{\pm};\qquad \widetilde H\mapsto H:=(a-bc)\widetilde H
\end{equation}
leads to the standard commutation relations, so we drop the tilde.
Set also $E=\sum u_i\del_i$.

Set $\text{weight}(u_1)=-\text{weight}(u_3)=w:=a-bc$,
$\text{weight}(u_2)=0$. Set
\begin{equation}\label{.19}
\fg_{-1}=\Tee(a, b, c)=\Span(\del_1, \del_2, \del_3)
\end{equation} and compute the Cartan prolong assuming
that $\fg_0$ is the smallest possible algebra containing $\fsl(2)$.

For $\un =(111)$, $\fg_2=0$; the space $\fg_1$ is easy to
get by hands; it is spanned by (for $a\neq 0$):
\scriptsize
\begin{equation}
\label{g10}
\renewcommand{\arraystretch}{1.4}
\begin{tabular}{|lc|c|}
\hline $\del_3^*:=$&$(a u_2^{(2)} + b c u_1 u_3)\del_1 +a(c u_1^{(2)}+ u_2
u_3)
\del_2 + (a c u_1 u_2 - (a + b c) u_3^{(2)}) \del_3$&$-w$\\
$\del_2^*:=$&$(u_1u_2+bu_3^{(2)})\del_1+(u_2^{(2)}+u_1u_3)\del_2+
 (u_2u_3+cu_1^{(2)})\del_3$&$0$\\
$\del_1^*:=$&$- ((a + b c) u_1^{(2)} -b u_2 u_3) \del_1
 + (a u_1 u_2 + b u_3^{(2)}) \del_2+ (b c u_1 u_3 + a~u_2^{(2)}) \del_3$&$w$\\
 \hline
\end{tabular}
\end{equation}
\normalsize
The commutators $[\fg_1, \fg_{-1}]$ are (here $E=1_3$)
\begin{equation}
\label{brakets}
\renewcommand{\arraystretch}{1.4}
\scriptsize{
\begin{tabular}{|c|c|c|c|}
\hline
&$\del_1^*$&$\del_2^*$&$\del_3^*$\\
\hline
$\del_1$&$H+aE$&$X^-$&$cX^+$\\
$\del_2$&$-X^+$&$E$&$aX^-$\\
$\del_3$&$bX^-$&$X^+$&$aE-H$\\
 \hline
\end{tabular}
}
\end{equation}
Since $[\del_2^*, \del_2]=E$, it follows that $\fg_{0}\simeq\fgl(2)$. Now, set $\fbr(2;
a, b, c):=\fg$.

\textbf{Occasional isomorphisms}. If $a\neq 0$, then we may assume that $a=1$ and then
\begin{equation}
\label{occiso}
\renewcommand{\arraystretch}{1.4}
\begin{array}{ll}
\fbr(2; 1, b, c)\simeq \fbr(2; 1, c, b),& \fbr(2; 1, 0, 0)\simeq
\fbr(2).
\end{array}
\end{equation}

\sssec{Particular cases of $\fbr(2; a, b, c)$}\label{ParCase} Since $\fbr(2; 1, 0, 0)= \fbr(2)$ has the
same non-positive part as $\fk(3; \One)$, it follows
that $\fbr(2)$ is a~partial Cartan prolong. All Kostrikin's examples
$L(\eps)$ also have the same non-positive parts as $\fk(3;
\One)$, so $L(\eps)$ is a~deformation $\fbr(2; a)$,
where $\eps=\frac{1+a}{a}$, 
of $\fo(5)$, see \cite{BLW}. Each
$L(\eps)$ can be embedded into $\fk(3; \One)$, see \cite{BLW}:
\begin{equation}\label{genb2}
\renewcommand{\arraystretch}{1.4}\scriptsize
\begin{array}{ll}
X_1^-=q^{(2)},&X_2^-=p,\\
X_1^+=-p^{(2)},&X_2^+=\delta(apq^{(2)} -qt) \text{ for
$\delta=\begin{cases}1&\text{if $a=2$,}\cr
\displaystyle\frac{1}{a+1}&\text{if $a\neq 2$.}\end{cases}$}\\
\end{array}
\end{equation}
Then,
\begin{equation}\label{.20}\scriptsize
h_1=pq,\quad h_2=\begin{cases}pq-t&\text{if $a=2$},\cr
\displaystyle\frac{a-1}{a+1}pq-\displaystyle\frac{1}{a+1}t&\text{if
$a\neq 2$}.\end{cases}
\end{equation}
The Cartan matrix is
\begin{equation}\label{.21}\scriptsize
\begin{pmatrix}
2&-1\cr
\alpha&0\cr
\end{pmatrix},\;\text{where $\alpha=\begin{cases}-1&\text{if $a=2$},\cr
\displaystyle\frac{a-1}{a+1}&\text{if $a\in \Kee\setminus\{1,
2\}$}.\end{cases}$}
\end{equation}
Kostrikin observed \cite{Ko} that $\fbr(2; a)$ can be deformed into
an algebra $L(\eps, \alpha, \beta)$ which may be identified with
$\fbr(2; a, b, c)$ for some $a, b, c$. For the complete description of
the deforms of $\fbr(2; \eps)$,
see \cite{BLW}.

\parbegin{Remark} There is no ``contact" analog of $\fbr(2; a, b,
c)$ because there is no non-degenerate antisymmetric bilinear form on
$\fg_{-1}=\Tee(a, b, c)$.
\end{Remark}

\sssbegin{Theorem}\label{ProloTabc} Let $\fg:=(\Tee(a, b, c), \fsl(2))_{*,\un}$ for the irreducible $\fsl(2)$-module
$\Tee(a, b, c)$, i.e., let ${a\neq bc}$. Then, $\fg_2\neq
0$ only for $a=0$. The parameter $c$ can be eliminated by setting $c=1$.
\end{Theorem}

Clearly, if $a=0$, we can divide $X_+$ by $b$, setting $b=1$. As
one can verify directly, $\fbr(2; 0, 1, c)_{\leq 0}$ has the
following $k$th Cartan prolong $\fg_k$ or rather $\fg_{k,\un}$, where $\un$ is not arbitrary but constrained: $\un =(11n)$, so we can set $\fg(c;\underline{n}):=(\Tee(0, 1, c), \fsl(2))_{*,\underline{n}}$. If $n=1$, then $\fg(c;\underline{1})_2=0$ whereas
 $\fg(c;\underline{1})_1\simeq (\fg(c;\underline{1})_{-1})^*$ is spanned by
\begin{equation}
\label{br01c}
\renewcommand{\arraystretch}{1.4}
\scriptsize{
\begin{tabular}{|lc|}
\hline $\del_1^*:=$&$(u_2u_3-cu_1^{(2)})\del_1 +
u_3^{(2)}\del_2 + c u_1 u_3 \del_3$\\
$\del_2^*:=$&$u_3^{(2)}\del_1$\\
$\del_3^*:=$&$- u_1u_3 \del_1+u_3^{(2)} \del_3$\\
 \hline
\end{tabular}
}
\end{equation}
The commutators $[\fg(c;\underline{1})_{-1}, \fg(c;\underline{1})_1]$ are given by the table
\begin{equation}
\label{braketsa=0}
\renewcommand{\arraystretch}{1.4}
\scriptsize{
\begin{tabular}{|c|c|c|c|}
\hline
&$\del_1^*$&$\del_2^*$&$\del_3^*$\\
\hline
$\del_1$&$h=c(-u_1\del_1+u_3 \del_3)$&$0$&$-x=-u_3 \del_1$\\
$\del_2$&$x=u_3 \del_1$&$0$&$0$\\
$\del_3$&$y=u_2\del_1+u_3\del_2+c u_1 \del_3$&$x=u_3 \del_1$&$\frac1ch=(-u_1\del_1+u_3 \del_3)$\\
 \hline
\end{tabular}
}
\end{equation}
Since $[\fg_1, \fg_{-1}]=\fsl(2)$, and $\fg_{\pm 1}$ are irreducible
$\fsl(2)$-modules, the  prolong ${\fbr(2; 0, 1,
c):=\mathop{\oplus}\limits_{i\geq -1}\fg_i}$ is a~simple Lie algebra.

It seems, $\fbr(2; (0, 1, c);\underline{(11n)})$ is a~new simple Lie algebra,
more precisely, it is a~deform of the non-standard Hamiltonian
algebra $\fh_{\exp}(2; \underline{(1, n)})$ preserving the form
$\omega = \exp(v)dv\wedge du$, and considered in \cite{BKK} in the
non-standard grading $\deg v=0$ and $\deg u=1$. The direct calculations, however, show that the parameter $c$ can be set equal to 1; so this is not a~new simple Lie algebra. (The situation can be illustrated with the following simplest example: although irreducible $\fgl(1)$-modules $V_\lambda$ corresponding to different characters $\lambda$ (i.e., the images of $1\in\fgl(1)$) are not isomorphic, the Cartan prolongations $(V_\lambda, \fgl(1))_*$ are all isomorphic to each other for $\lambda \neq 0$.)

\ssec{Skryabin algebras} To describe them,
observe that, as far as $\fvect(3, \un)$-action is
concerned, we have the following identifications: $v^2=v^{-1}$ and
(cf. (\ref{ident1}))
\begin{equation}\label{dudu/v3}
du_idu_jv^{-1}=\sign(ijk)\del_k\;\text{ for any permutation $(ijk)$
of $(123)$}.
\end{equation}
Observe that, as $\fvect(3; \un)$-modules, $\fsvect(3;
\un)\not\simeq Z^2(3; \un)$, more precisely,
$\fsvect(3, \un)$ is not a~$\fvect(3;
\un)$-module: being a~simple Lie algebra, $\fvect(3;
\un)$ has no submodules in the adjoint representation.
Observe that neither $\fsvect(3;
\un)$ nor $Z^2(3; \un)$ are $\cO(3;
\un)$-modules.

\sssec{The deep Skryabin algebra} As $\Zee/4$-graded
Lie algebras, we have, see \cite{Sk3}:
\begin{equation}\label{defdy}\scriptsize
\renewcommand{\arraystretch}{1.4}
\begin{array}{l}
\fdy(\un):=\fg_\ev \oplus\fg_{\bar
1}\oplus\fg_{\bar 2}\oplus\fg_{\bar 3}
\simeq \\
\fvect(3; \un)\oplus \cO(3; \un)_{-\Div}\oplus
\Omega^1(3;\un)_{\Div}\oplus Z^2(3;
\un) \ ,\\
\fdy(\un)^{(1)} \simeq \fvect(3; \un)\oplus
\cO(3; \un)_{-\Div}\oplus
\Omega^1(3;\un)_{\Div}\oplus B^2(3; \un) \ .
\end{array}
\end{equation}
In particular (hereafter $|\un |=\sum N_i$),
\begin{equation}\label{dimdy}
\dim \fdy(\un)= 3^{|\un |+2}+1,\quad \dim
\fdy(\un)^{(1)}= 3^{|\un |+2}-2.
\end{equation}

\textit{The following questions should be applied to other Skryabin algebras as well,
especially for $\un=\One$: ``Which one is simple: the algebra or only its derived algebra?
Is this indeed the first derived algebra that is simple, or we have to pass
to the second derived algebra?''} 

The multiplication in $\fdy(\un)$ is given
by the $\fvect(3, \un)$-invariant bilinear differential operators
acting in the spaces of tensor fields entering (\ref{defdy}). Over
$\Cee$, all such operators are described \cite{Gr0}; to describe
even unary operators is an open problem for $p>0$. For $p=3$, the
following formulas for multiplication reveal presence of new
invariant operators; for more examples, see \cite{BjL1}.

In eq. \eqref{multdy}, we have $f, g\in\cO$ and $\omega^i\in\Omega^i$; in lines
3 and 5 we postulate antisymmetry of the brackets, cf.
(\ref{multme}):
\begin{equation}\label{multdy}\scriptsize
\renewcommand{\arraystretch}{1.4}
\begin{array}{lcl}
{}[fv^{-1}, gv^{-1}]&=& (gdf-fdg)v;\\
{}[fv^{-1}, \omega^1v]&=&
-d(f\omega^1);\\
{}[fv^{-1}, \omega^2]&=&f\omega^2v^{-1}\in\fg_\ev;\\
{}[\omega^1_1v, \omega^1_2v]&=&
\omega^1_1\omega^1_2v^{-1}\in\fg_\ev;\\
{}[\omega^1v, \omega^2]&=&
\omega^1\omega^2v\in\cO_{2\Div}\simeq\cO_{-\Div};\\
{}[\sum f_i(u)du_jdu_k, \sum g_i(u)du_jdu_k]&=&
(f_3g_2-f_2g_3)du_1v + cycle(123).
\end{array}
\end{equation}
Set
(\cite{Sk3}):
\begin{equation}\label{grdy}
\renewcommand{\arraystretch}{1.4}
\begin{array}{lll}
\deg u^{(\underline{r})}\del_i =4|\underline{r}|-4,&&
\deg u^{(\underline{r})}du_iv =4|\underline{r}|-2,\\
\deg u^{(\underline{r})}v^{-1} =4|\underline{r}|-3,&&
\deg u^{(\underline{r})}du_idu_j =4|\underline{r}|-1.\\
\end{array}
\end{equation}
For $\fgl(3)=\fgl(V)$, where $V:= V[-1]$, we have 
\[
\text{$E^2(V)\simeq
(V[-1])^*\simeq V^*[-2]$, and $E^3(V)\simeq \One[-3]$,}
\]
 where $\One$
is the trivial $\fsl(3)$-module. Then, $\fdy(\un)$ can be defined as
the Cartan-Tanaka-Shchepochkina prolong with the following
non-positive part (here $i, j=1, 2, 3$):
\begin{equation}\label{DY}
\renewcommand{\arraystretch}{1.4}
\scriptsize{
\begin{tabular}{|c|c|c|c|c|c|}
\hline &$\fg_0$&$\fg_{-1}$&$\fg_{-2}$&$\fg_{-3}$&$\fg_{-4}$\cr
\hline
$\Span_\Kee$&$x_j\del_i$&$du_idu_j$&$v\,du_i$&$v^{-1}$&$\del_i$\cr
\hline &$\fgl(3)\simeq
\fgl(V)$&$V=V[-1]$&$E^2(V)$&$E^3(V)$&$V[-4]$\cr \hline
\end{tabular}
}
\end{equation}

\sssec{The big Skryabin algebra} This Lie algebra, $\fby(\un)$, was
only described so far (\cite{Sk3, S}) as having the following
non-positive part:
\begin{equation}\label{BY}
\renewcommand{\arraystretch}{1.4}
\scriptsize{
\begin{tabular}{|c|c|c|c|}
\hline $\fg_0$&$\fg_{-1}$&$\fg_{-2}$&$\fg_{-3}$\cr
\hline
$\fgl(3):=\fgl(V)$&$V$&$E^2(V)$&$E^3(V)$\cr
\hline
\end{tabular}
}
\end{equation}

\sssec{The middle Skryabin algebra} As $\Zee/2$-graded
Lie algebras, we have
\begin{equation}\label{defmy}
\fmy(\un)\simeq \fvect(3; \un)\oplus \Omega^1(3;
\un)_{\Div}.
\end{equation}
with multiplication given by (\ref{multdy}).
Set (\cite{S}):
\begin{equation}\label{grmy}
\renewcommand{\arraystretch}{1.4}
\begin{array}{l}
\deg u^{(\underline{r})}\del_i =2|\underline{r}|-2,\ \ \ \
\deg u^{(\underline{r})}du_iv =2|\underline{r}|-1.\\
\end{array}
\end{equation}
Then, $\fmy(\un)$ can be defined as the generalized prolong with the
following non-positive part (here $i, j=1, 2, 3$):
\begin{equation}\label{MY}
\renewcommand{\arraystretch}{1.4}
\scriptsize{
\begin{tabular}{|c|c|c|c|}
\hline &$\fg_0$&$\fg_{-1}$&$\fg_{-2}$\cr \hline
$\Span_\Kee$&$u_j\del_i$&$du_iv$&$\del_i$\cr \hline
&$\fgl(3)\simeq\fgl(V)$&$V$&$E^2(V)$\cr \hline
\end{tabular}
}
\end{equation}
and therefore the non-positive part $\fmy(\un)_{\leq 0}$ coincides
with $\fo(7)_{\leq 0}$ in the grading with the last ``selected\rq\rq
root.

\sssec{The little Skryabin algebras} The dimensions of these algebras
$\fg$ are distinct, as well as the structures they preserve, so they
should be considered as separate entities. Here we only consider one
of these algebras --- $\text{Y}_{(1)}(\un)$ in the original
notation --- defined as the generalized prolong with the following
non-positive part:
\begin{equation}\label{my''}
\renewcommand{\arraystretch}{1.4}
\scriptsize{
\begin{tabular}{|c|c|c|c|} \hline
&$\fg_0$&$\fg_{-1}$&$\fg_{-2}$\cr \hline $\Span_\Kee$&$\fsl(3)\simeq
\fsl(V)$& $V=\Span(du_i)_{i\leq 3}$&$E^2(V)=\Span(\del_i)_{i\leq
3}$\cr \hline
\end{tabular}
}
\end{equation} Clearly, $Y_{(1)}(\un)\simeq
\fsmy(\un)$ and, as $\Zee/2$-graded Lie algebras,
\begin{equation}\label{defsmy}
\renewcommand{\arraystretch}{1.4}\scriptsize
\begin{array}{l}
\fsmy(\un)\simeq \fsvect(3; \un)\oplus
Z^1(3; \un),\\
\fsmy^{(1)}(\un)\simeq \fsvect(3;
\un)^{(1)}\oplus B^1(3; \un),
\end{array}
\end{equation}
so
\begin{equation}\label{dimSMY}
\dim \fsmy^{(1)}(\un)=3^{|\un
|+1}-3.
\end{equation}

There are known 5 little Skryabin algebras (\cite{Sk3, S}); the other
little algebras are not graded but filtered deforms studied by Ladilova and further considered in \cite{GKL}.

\ssec{Frank algebras $\ffr(n)$}\label{FrGen} The non-positive part of $\ffr(n)$ is
the same as that of $\fg:=\fk(3; (n,1, 1))$, where the coordinates
of the vector $\un =(n,1, 1)$ correspond to the ordered set $(t;q,
p)$; and hence same non-positive part as $\fsp(4)_{\leq 0}$ in the
grading with the first ``selected\rq\rq root. The bases of
components of $\ffr(n)_i$ for $i>0$ are as follows:
\begin{equation}\label{.40}\scriptsize
\begin{array}{ll}
\text{for }\fg_{2i-1}:&\left\{p{}^{(2)}q t^{(i-1)}-pt^{(i)},\ \
q{}^{(2)}p t^{(i-1)}+q t^{(i)}\right\},\\

\text{for }\fg_{2i}:&\left\{p{}^{(2)}t^{(i)}, \ \
pqt^{(i)}, \ \ q{}^{(2)}t^{(i)},\ \ p{}^{(2)}q{}^{(2)}t^{(i-1)}-t^{(i+1)} \right\}
\text{~~if $1\leq i\leq \mathcal{N}$}\\

\text{for }\fg_{2\mathcal{N}+2}:&\left\{p{}^{(2)}t^{(\mathcal{N}+1)}, \ \
pqt^{(\mathcal{N}+1)}, \ \ q{}^{(2)}t^{(\mathcal{N}+1)} \right\},\\\end{array}
\end{equation}
where the first line is for $1\leq i\leq \mathcal{N}+1$ and $\mathcal{N}:=3^n-2$, so that $t^{(\mathcal{N}+1)}$ is the top power of $t$ allowed by $\un=(n,1,1)$; hence
\[
\dim\ffr(n)=2\cdot 3^{n+1}.
\]
Here the functions generate contact vector fields by the bracket with a~plus in front of the Poisson bracket:
\be\label{kbfr}
\{f,g\}=(2-E)(f)\del_t(g)-\del_t(f)(2-E)(g)+\del_{p}(f)\del_{q}(g)- \del_{q}(f)\del_{p}(g),
\ee
the one corresponding to the form $dt-p\,dq+q\,dp$.

\normalsize In particular, as $\fg_0(=\fsp(2)\oplus \Kee z)$-module, $\fk(3; \One)_1$ has two lowest
weight vectors: $qt$ and $pq^{(2)}$, whereas
\begin{equation}\label{.22}
\ffr(n)_1=\Span(tp-p^{(2)}q,\ tq+pq^{(2)})\subset \fk(3; \One)_1.
\end{equation}

\sssbegin{Remark}\label{remFr} For the majority of simple
Weisfeiler graded (see
\cite{LSh}) vectorial Lie (super)al\-geb\-ras~ $\fg$, their
positive part is (at least for $\un=\One$ and $p>3$, or $p=0$)
generated by $\fg_1$. Among the simple $\Zee$-graded Lie superalgebras in characteristic 0 or $p>3$ described in \cite{BGLLS2}, exceptions whose positive part is not generated by $\fg_1$ are $\fvect(1;\underline{1}|0)$,
$\fk(1;\underline{1}|1)$,  and $\fk(1;\underline{1}|2)$. 

The Frank algebras are also such
exceptions: the non-positive components of $\ffr(n)$ and $\ffr_1$ for
any $n$ generate $\fsp(4)$. To generate $\ffr(n)$, we should add the following
lowest weight operators (with respect to $\ffr(n)_0=\fgl(2)$): \scriptsize
\begin{equation}\label{.fr40}
\begin{array}{ll}
\text{$z_1= pq^{(2)} + qt$ and $z_2= q^{(2)}t$}&\text{for $n=1$},\\

\text{$z_1, \dots, z_{i+2}=t^{(3^i)} - p^{(2)}q^{(2)}t^{(3^i-2)}$}&\text{for $1\leq i< n$}.\\
\end{array}
\end{equation}
\normalsize
\end{Remark}

\ssec{Ermolaev algebra $\fer$}\label{Er} As $\Zee/2$-graded Lie
algebra, it is defined as follows:
\begin{equation}
\label{er1}\scriptsize
\renewcommand{\arraystretch}{1.4}
\begin{array}{l}
\fer(\un):=\fvect(2; \un)\oplus
\cO(2;\un)_{\Div},\\
\fer(\un)^{(1)}= \fvect(2; \un)\oplus
\cO'(2;\un)_{\Div}, \text{~~so $\dim\fer(\un)^{(1)}=3^{|\un |+1}-1$.}
\end{array}\end{equation}

To define the bracket, recall that $v^2=v^{-1}$ and observe that (cf. (\ref{dudu/v3}))
\begin{equation}\label{dudu/v2}
du_iv^{-1}=\sign(ij)\del_j\;\text{ for any permutation $(ij)$ of
$(12)$}.
\end{equation}
For any $f\ v, g\ v \in\cO(2; \un)_{\Div}$, set
\begin{equation}
\label{er2} {}[fv, gv]:=(fdg-gdf)v^{-1}= (f\del_2(g)-g\del
_2(f))\del_1+(g\del_1(f)-f\del _1(g))\del_2 \text{};
\end{equation}
define the other products canon-ically. Define the $\Zee$-grading of
$\fer(\un)$ as induced by the standard $\Zee$-grading of
$\cO(2;\un)$:
\begin{equation}
\label{er3} \fer(\un)_i:=\fvect(2; \un)_i\oplus
\left(\cO(2;\un)_{\Div}\right)_{i+1}\quad \text{for any
$i\geq -1$}.
\end{equation}
Then, $\fer(\widetilde\un)$ is defined as the Cartan prolong of the
following non-positive part:
\begin{equation}
\label{er}\scriptsize
\renewcommand{\arraystretch}{1.4}
\begin{array}{lll}
\fer_{-1}&=&\Span_\Kee(\del_1, \del_2;\; \vvol);\\
\fer_{0}&=&\Span_\Kee(u_i\del_j;\; u_k\vvol\mid i, j, k=1, 2).
\end{array}
\end{equation}


\ssec{Questions we address for all the Lie algebras described in this section}\label{Qs} In this paper, we managed to answer only some of them.

1) What are
the structures these algebras preserve?

2) What are the complete and partial Cartan-Tanaka-Shchepochkina
prolongs of the non-positive and negative parts of $\fg$
corresponding to the $\Zee$-gradings obtained by setting $\deg
X_i^{\pm}=\pm \delta_{i, i_0}$ for one or several selected indices
$i_0$?

3) What are the defining relations explicitly?

4) What are the ``most natural'' generators of $\fbr(2; a, b, c)$, and the corresponding relations? Conjecturally, the answer is similar to
\cite{GL5}.

\vskip 0.2 cm

The above problems are resolved in what follows, at least, partly.
Questions 5)--7) are open:

\vskip 0.2 cm

5) The Brown algebras $\fg=\fbr(r)$
given by tridiagonal $r\times r$ Cartan matrices
of type $1\fbr(3)$, see (\ref{br3a}), of size $r>3$ with one 0 on the main diagonal seem to be of infinite-dimension, though $\dim \mathop{\oplus}\limits_{|i|\leq n}\fg_i$, where the $\Zee$-grading
of these Lie algebras of the form $\fg(A)$ is defined by setting $\deg X_i^{\pm}=\pm 1$ for
all the
Chevalley generators $X_i^{\pm}$, grows rather slow as
$n\tto \infty$. Which of these
$\Zee$-graded algebras are of polynomial growth?
\textbf{Conjecturally, these algebras
grow polynomially in $n$ for $r=4, 5$}.

6) What are the analogs of Weisfeiler gradings (see
\cite{LSh}) for infinite-dimensional simple vectorial Lie algebras in the limit as variable (non--critical)
coordinates of $\un $ tend to $\infty$? (The corresponding
preserved structures are most interesting.)

7a) In what follows we impose no restrictions on $\un $; the
construction of prolongation imposes them automatically.
Constructing the prolongs of $\mathop{\oplus}\limits_{-d\leq i\leq
0}\fg_i$ when $d>1$, we observed that if $\un _j$ can be $>1$, then $\deg
u_j>1$. What is the reason for this? The opposite is not a~must, as
the Melikyan algebras for $p=5$ show, see \eqref{.23} and Subsection~
\ref{ConcealMel}.

We thought that if $d=1$, then either $N_i=1$ for all $i$, or no
restrictions on $N_i$ for all $i$ ``due to interchangeability of
indeterminates for the case of transitive (in particular, simple)
vectorial Lie algebras". There are, however, examples, where we observe symmetry, in other words ``interchangeability" of some
indeterminates, whereas some of the coordinates of $\un$ corresponding to these indeterminates are constrained
(critical) whereas several other ones are not for no reason we can formulate.

7b) How to find a~concealed parameter? Observe that from the
realization \eqref{er1} one may deduce that the shearing vector depends on 2
parameters, not on 3 parameters as decomposition \eqref{er} hints. In Subsection~ \ref{concER} we will see that, indeed, $\Par\un=2$ for $\fer(\un)$.

It was only thanks to a~computer experiment that we have discovered
the correct number of parameters the shearing vector $\un$ for $\fby(\un)$ depends on.
This discovery should not have surprised us since all indeterminates of degree 2,
see \eqref{realby}, are on equal footing.

8) What is the complete list of deforms of the ``standard" simple Lie algebras and their ``relatives",
cf. \cite{BLLS1}?

\section{Interpretations}

\ssec{Melikyan algebras for $p=5$} We realize $\fg(2)_-\oplus
\fg(2)_0$ by vector fields, thanks to Shchepochkina's algorithm
\cite{Shch}, as follows: \scriptsize
\begin{equation}\label{realme}
\renewcommand{\arraystretch}{1.4}
\begin{tabular}{|l|l|}
\hline $\fg_{0}$&$X_- = u_2 \del_1 + u_1 u_2^{(2)} \del_4 + u_2^{(3)} \del_5
+ u_5
\del_4$,\\
&$X_+ = u_1 \del_2 + 2 u_1^{(3)} \del_4 + u_4 \del_5$,\\
&$h_1 = u_1 \del_1 + u_3 \del_3 + 2 u_4 \del_4 + u_5 \del_5$,\\
&$h_2 = u_2 \del_2 + u_3 \del_3 + u_4 \del_4 + 2 u_5 \del_5$\\
\hline\end{tabular}\quad \renewcommand{\arraystretch}{1.4}
\begin{tabular}{|l|l|}
\hline $\fg_{-1}$&$\del_1 -u_2 \del_3- u_1 u_2 \del_4 - u_2^{(2)}
\del_5$;\quad
$\del_2$\\
\hline
$\fg_{-2}$&$\del_3 + u_1 \del_4 + u_2 \del_5$\\
\hline
$\fg_{-3}$&$\del_4$;\quad $\del_5$\\
\hline
\end{tabular}
\end{equation}
\normalsize Set $\deg u_1=(1, 0)$, $\deg u_2=(0, 1)$. This
determines the other degrees: $\deg u_3=(1, 1)$, ${\deg u_4=(2, 1)}$,
$\deg u_5=(1, 2)$. Unlike $p=0$ case, the complete prolong of
$\fg(2)_-\oplus \fg(2)_0$ in the depth 3 grading of $\fg(2)$
strictly contains $\fg(2)$ (the underlined components) and has (for
the simplest $\un $) the following irreducible components
as $\fg_0$-modules given by their highest weights: \tiny
\begin{equation}\label{.23}
\renewcommand{\arraystretch}{1.4}
\begin{tabular}{|l|l|l|l|l|l|}
\hline $\deg$&$\dim$&highest weights&$\deg$&$\dim$&highest weights
\cr \hline $-3$&$2$&$\underline{(-1,-2)}$&$23$&$2$&$(12, 11)$\cr
$-2$&$1$&$\underline{(-1,-1)}$&$22$&$1$&$(11, 11)$\cr
$-1$&$2$&$\underline{(0,-1)}$&$21$&$2$&$(11, 10)$\cr
$0$&$4$&$\underline{(1,-1)}, \underline{(0,0)}$&$20$&$4$&$(11, 9),
(10, 10)$\cr $1$&$2$&$\underline{(1,0)}$&$19$&$2$&$(10,9)$\cr
$2$&$4$&$(2,0), \underline{(1,1)}$&$18$&$4$&$(10,8)$, $(9,9)$\cr
$3$&$6$&$(3,0)$, $\underline{(2,1)}$&$17$&$6$&$(10, 7), (9, 8)$\cr
$4$&$3$&$(3,1)$&$16$&$3$&$(9,7)$\cr $5$&$6$&$(4,1)$,
$(3,2)$&$15$&$6$&$(9,6), (8,3)$\cr
$6$&$8$&$(5,1)$,$(4,2)$&$14$&$8$&$(9,5)$, $(8,7)$\cr
$7$&$4$&$(5,2)$&$13$&$4$&$(8,5)$\cr $8$&$8$&$(6, 2)$,
$(5,3)$&$12$&$8$&$(8,4), (7,5)$\cr $9$&$10$&$(7,2)$,
$(6,3)$&$11$&$10$&$(8,3)$, $(7, 4)$\cr $10$&$5$&$(7,3)$&&&\cr \hline
\end{tabular}
\end{equation}
\normalsize

\sssec{On concealed parameters of $\protect\un$}\label{ConcealMel} Computer experiments show that
without restrictions on
$\un $ the complete prolong only depends on two
parameters, as predicted by eqs.~\eqref{defme}, i.e., $\un =(1,1,1,
N_4, N_5)$. Verdict: no concealed parameters.

\sssec{Melikyan algebras for $p<5$ and their superizations} See \cite{BjL,
BGL, BGLLS}.

\ssec{Brown algebras} Let $x^{\pm}_i$ be the preimage of the
generators $X^{\pm}_i$ relative (\ref{myst}). By abuse of notation
we will often write $x^{\pm}_i$ instead of $X^{\pm}_i$; let $x_i$ be
either all $x^{+}_i$ or all $x^{-}_i$.

\underline{$\fbr(2)$}: Basis (of $\fbr(2)_{\pm}$):
\begin{equation}\label{.24}
x_1,\; x_2,\; [x_1,x_2],\; [x_2,[x_2, x_1]].
\end{equation}
Here, $x^{+}_2$ and $x^{-}_2$ generate $\fhei(2; 3; \underline{1})$
on which $h_1$ acts as an outer derivation. The Fock space
representation $\cO(1; \underline{1})$ of $\fhei(2; 3;
\underline{1})\inplus \Kee h_1$ (hereafter $\fa\niplus\fii$ is a
semidirect sum of algebras, where $\fii$ is an ideal) is irreducible
of dimension 3. Therefore, the non-positive terms of the simplest
$\Zee$-gradings ($\deg x^{\pm}_{i_0}=\pm 1$) are:
\begin{equation}\label{b2}
\renewcommand{\arraystretch}{1.4}
\scriptsize{
\begin{tabular}{|c|c|c|c|}
\hline
$i_0$&$\fg_0$&$\fg_{-1}$&$\fg_{-2}$\cr
\hline
$1$&$\fhei(2;
3; \underline{1})\inplus \Kee h_1$&$\cO(1; \underline{1})$&$-$\cr
\hline
$2$&$\fgl(2)\simeq\fgl(V)$&$V$&$E^2(V)$\cr
\hline
\end{tabular}
}
\end{equation}

\ssec{Remarks} The first grading tempts us to investigate if there is a~non-trivial
Cartan prolong of the pair $\fg_0=\fhei(2; 3; \un)\inplus
\Kee h_1$ and $\fg_{-1}=\cO(1; \un)$. But for the prolong
to be simple, irreducibility is needed, while $\cO(1;
\un)$ is irreducible $\fhei(2; 3; \un)$-module
only for $N=1$.

Observe that in the second grading, the non-positive terms are the
same as the non-positive terms of $\fk(3)$, and hence same as those
of the Frank algebras and same as those of $\fsp(4)$ with one
(first) ``selected\rq\rq simple root. Moreover, even the defining
relations between the positive (negative) generators are the same as
the Serre relations of $\fsp(4)$ although the Cartan matrix of
$\fbr(2)$, to say nothing of $\fbr(2;\eps)$, is different.

\underline{$\fbr(3)$}: The Chevalley bases of $1\fbr(3)_{\pm}$ and
$2\fbr(3)_{\pm}$ are, respectively, as follows:
\begin{equation}\label{genb3}
\tiny
\begin{array}{llll}
\renewcommand{\arraystretch}{1.4}
1)\begin{pmatrix}2&-1&0\\-1&2&-1\\
0&-1&0\end{pmatrix}:&\begin{array}{l}
x_1,\; x_2,\; x_3;\\
{}[x_1,x_2],\quad [x_2,x_3];\\
{}[x_3,[x_3, x_2]],\quad [x_3,[x_2,x_1]];\\
{}[x_3,[x_3,[x_1,x_2]]];\\
{}[[x_2,x_3], [x_3,[x_1,x_2]]];\\
{}[[x_3,[x_1,x_2]], [x_3,[x_2,x_3]]];\\
{}[[x_3,[x_2,x_3]], [x_3,[x_3,[x_1,x_2]]]];\\
{}[[x_3,[x_2,x_3]], [[x_2,x_3], [x_3,[x_1,x_2]]]];\\
{}[[x_3,[x_3,[x_1,x_2]]], [[x_2,x_3], [x_3,[x_1,x_2]]]].\end{array}
&2)\begin{pmatrix}2&-1&0\\-2&2&-1\\
0&-1&0\end{pmatrix}:&
\renewcommand{\arraystretch}{1.4}
\begin{array}{l}
x_1,\; x_2,\; x_3;\\
{}[x_1, x_2],\ [x_2,x_3],\\
{}[x_2, [x_1, x_2]],\ [x_3, [x_1, x_2]],\ {}[x_3, [x_2, x_3]],\\
{}[x_3, [x_3, [x_1, x_2]]],\
{}[[x_1, x_2], [x_2, x_3]],\\
{}[[x_2, x_3], [x_3, [x_1,x_2]]],\\
{}[[x_2, [x_1, x_2]], [x_3, [x_2, x_3]]],\\
{}[[x_3, [x_1, x_2]], [[x_1, x_2], [x_2, x_3]]].\end{array}
\end{array}
\normalsize
\end{equation}
For $\fg=1\fbr(3)_{\pm}$, the non-positive terms of the $\Zee$-gradings
in terms of $\fg_0$-modules are as follows (underlined are the
dimensions of the irreducible $\fg_0$- or $\fa$- and $\fb$-modules if $\fg_0=\fa\niplus\fb$): \begin{equation}\label{.25}
\renewcommand{\arraystretch}{1.4}
\scriptsize{
\begin{tabular}{|c|c|c|c|c|c|}
\hline
$i_0$&$\fg_0$&$\fg_{-1}$&$\fg_{-2}$&$\fg_{-3}$&$\fg_{-4}$\cr
\hline
$1$&$\fbr(2)\inplus \Kee
h_1$&$\underline{8}$&$\underline{1}$&$-$&$-$\cr
\hline
$2$&$\left(\fsl(2)\inplus \Kee h_2\right)\niplus\fhei(2; 3;
\underline{1})$&$\underline{2}\boxtimes
\underline{3}$&$\underline{1}\boxtimes\underline{3}$&$\underline{2}\boxtimes\underline{1}$
&$-$\cr
\hline
$3$&$\fgl(3)$&$\underline{3}$&$\underline{3}$&$\underline{1}$&$\underline{3}$\cr
\hline
\end{tabular}
}
\end{equation}
\textbf{The last line coincides with $\fdy_-$, see (\ref{DY}); the
first line shows that $1\fbr(3)$ is a~partial prolong of
$(\fsp(10)_-, \fbr(2)\inplus \Kee h_1)$ in the contact grading of
$\fsp(10)$. The prolongs of the non-positive parts in lines $1$ and $2$ return $1\fbr(3)$.}

For $\fg=2\fbr(3)_{\pm}$, the non-positive terms of the $\Zee$-gradings
in terms of $\fg_0$-modules are as follows (to be used later:
$\fg_{-1}$ has HWV of weight $(0,1,1)$ and $\fg_{-2}$ has a~HWV of
weight $(1,0,2)$.) The prolongs of the non-positive parts in all three lines return
$2\fbr(3)$: the Cartan prolong of $\fg_-\oplus\fg_0$ has $\fg_1$ and $\fg_2$ of the
dimensions they have in $\fbr(3)$ and $\fg_3=0$. In line $3$, this is
$\fg_1=\underline{4}$ and $\fg_2=\underline{5}$, so $\fbr(3)$ is the Cartan prolong of
$(\fg_-,\fc\fsp(4))$, where $\fg_-=\fg_{-1}\oplus\fg_{-2}$ is of dimension $9$ with
$\fg_{-2}=[\fg_{-1},\fg_{-1}]$. 
\begin{equation}\label{.251}
\renewcommand{\arraystretch}{1.4}
\scriptsize{
\begin{tabular}{|c|c|c|c|c|}
\hline $i_0$&$\fg_0$&$\fg_{-1}$&$\fg_{-2}$&$\fg_{-3}$\cr \hline
$1$&$\fbr(2)\inplus \Kee(h_3-
h_1)$&$\underline{8}$&$\underline{1}$ of weight $(-1,0,1)$&$-$\cr \hline
$2$&$\left(\fsl(2)\inplus \Kee h_2\right)\niplus\fhei(2; 3;
\underline{1})$&$\underline{2}\boxtimes
\underline{3}$&$\underline{1}\boxtimes\underline{3}$&$\underline{2}\boxtimes\underline{1}$\cr
\hline
$3$&$\fc\fsp(4)\simeq\fc\fo(5)$&$\underline{4}$&$\underline{5}$&$-$\cr
\hline
\end{tabular}
}
\end{equation}

\ssec{The deep Skryabin algebra $\fg=\fdy$}\label{fdy} Let $\fg_-$ be realized by vector fields as follows:
\scriptsize
\begin{equation}\label{realdy2}
\renewcommand{\arraystretch}{1.4}
\begin{tabular}{|l|l|}
 \hline $\fg_{-1}$
 &$X_1 = 2 \del_1 + 2u_2\del_4 + 2u_3\del_5 + 2u_2^{(2)}u_3\del_9 + (u_2u_3 + 2u_6)\del_7 + (2u_3u_4 + u_7)\del_8$,\\
 &$X_2 = 2 \del_2 + u_5\del_7 + u_7\del_9,\qquad
 X_3 = 2 \del_3 + u_2\del_6 + 2u_4\del_7 + u_2u_4\del_9 + (u_2 u_5+u_7)\del_{10}$\\
 \hline $\fg_{-2}$
 &$X_4:= \del_{4} + u_{5} \del_{8} + u_{6} \del_{9}$,\qquad
 $X_5:= \del_{5} + u_{6}\del_{10}$,\qquad
 $X_6:= \del_{6}$\\
 \hline
 $\fg_{-3}$
 &$X_7:=\del_{7}$\\
 \hline $\fg_{-4}$
 &$X_8:=\del_{8},\quad X_9:=\del_{9},\quad
 X_{10}:=\del_{10}$\\
\hline
\end{tabular}
\end{equation}
\normalsize
Then, $\fdy_0=\fgl(3)$ and $\fs\fdy_0=\fsl(3)$ are spanned by
\scriptsize
\begin{equation}\label{realdy}
\renewcommand{\arraystretch}{1.4}
\begin{tabular}{|ll|}
\hline
$x_1 =$&$ u_3\del_2 + 2u_3^{(2)}\del_6 + (2u_4u_3^{(2)} + u_{10})\del_9 + u_5\del_4 + 2u_3^{(2)}u_5\del_{10} + u_5^{(2)}\del_8$,\\
$x_2 =$&$ 2u_2\del_1 + 2u_2^{(2)}\del_4 + 2u_2^{(2)}u_3\del_7 + (2u_2u_3 + 2u_6)\del_5 + (2u_2u_3u_4 + 2u_9)\del_8 + 2u_6^{(2)}\del_{10}$,\\
$x_3 =$&$ [x_1,x_2] = 2u_3\del_1 + 2u_3^{(2)}\del_5 + u_2u_3^{(2)}\del_7 + (2u_2^{(2)}u_3^{(2)} + u_6^{(2)})\del_9 + (2u_4u_3^{(2)} + u_5u_6 + 2u_{10})\del_8 + u_6\del_4$,\\
$y_1 =$&$ u_2\del_3 + 2u_2^{(2)}\del_6 + u_4\del_5 + 2u_2^{(2)}u_4\del_9 + u_4^{(2)}\del_8 + (2u_5u_2^{(2)} + u_9)\del_{10}$,\\
$y_2 =$&$ 2u_1\del_2 + 2u_1^{(2)}\del_4 + 2u_1^{(2)}u_3\del_7 + (u_1u_3 + 2u_5)\del_6 + (u_2u_3u_1^{(2)} + u_3u_4u_1 + 2u_8)\del_9$\\
 &$ + (u_1^{(2)}u_3^{(2)} + u_1u_5u_3 + 2u_5^{(2)})\del_{10}$,\\
$y_3 =$&$ [y_1,y_2] = u_1\del_3 + u_1^{(2)}\del_5 + 2u_1^{(2)}u_2\del_7 + (u_1^{(2)}u_2^{(2)} + 2u_4^{(2)})\del_9 + 2u_4\del_6
 + u_1^{(2)}u_4\del_8 + (2 u_4 u_5+u_8)\del_{10}$,\\
$h_1 =$&$ [x_1,y_1] = 2u_2\del_2 + u_3\del_3 + 2u_4\del_4 + u_5\del_5 + 2u_9\del_9 + u_{10}\del_{10}$,\\
$h_2 =$&$ [x_2,y_2] = 2u_1\del_1 + u_2\del_2 + 2u_5\del_5 + u_6\del_6 + 2u_8\del_8 + u_9\del_9$,\\
$h_3 =$&$2u_2\del_2 + 2u_3\del_3 + 2u_4\del_4 + 2u_5\del_5 + u_6\del_6+u_7\del_7+u_8\del_8$\\
\hline
\end{tabular}
\end{equation}
\normalsize
The weights in \eqref{.26} and \eqref{1.1} are given with respect to the following diagonal elements of~$\fdy_0$:
\[\scriptsize
\begin{array}{l}
H_1:=u_1\del_1+u_4\del_4+u_5\del_5+u_7\del_7-u_8\del_8+u_9\del_9+u_{10}\del_{10},\\
H_2:=u_2\del_2+u_4\del_4+u_6\del_6+u_7\del_7+u_8\del_8-u_9\del_9+u_{10}\del_{10},\\
H_3:=u_3\del_3+u_5\del_5+u_6\del_6+u_7\del_7+u_8\del_8+u_9\del_9-u_{10}\del_{10}.
\end{array}
\]

Let us describe $\fdy$ in the form similar to $\fme(\un)$, i.e., as a~sum
of $\fvect(3, \un)$ and modules over it. The lines $\deg=-4$ through
$-1$ in table \eqref{1.1} give us weights with respect to the $h_i$,
see \eqref{realdy}. It is clear that $\deg=-3$ corresponds to the
volume forms and either line $-4$ or line $-2$ should correspond to
$\fvect(3, \un)$ which should lie in even degrees. A few
experiments approve just one scenario. Line $-4$ gives us the
highest weight of $\fvect_{-1}(3, \un)$ in a~non-standard grading:
$\text{weight}(\del_{3}) = (-1,-1, -2)$ (which implies that
$\text{weight}(x_1) = (2,1,1)$, and hence $
\text{weight}(\vvol)=(1,1,1)$), and hence the highest weights and
the corresponding vectors of the following components are as
follows: \scriptsize
\begin{equation}\label{.26}
\renewcommand{\arraystretch}{1.4}
\begin{array}{ll}
\fg_{-1}:& \text{weight}(du_1du_2) = (3,3,2) \cong(0,0,-1),\\
\fg_{-2}: &\text{weight}(du_1\cdot\vvol) =(3,2,2) \cong (0,-1,-1),\\
\fg_{-3}:& \text{weight}(\vvol^{-1}) = (-1,-1,-1),\\
\fg_{-4}:& \text{weight}(\del_3) = (-1,-1,-2).
\end{array}
\end{equation}
\normalsize
\sssec{A concealed parameter}\label{concDY} The  realization of
$\fdy(\un)$, see \eqref{1.1}, shows that $\un \in\Nee^{10}$. Representation
(\ref{defdy}) shows that $\un \in\Nee^{10}$ depends on at least 3
parameters. The concealed parameter we found for $\fby(\un)$ urges us
to make sure we have established $\Par\un$ for $\fdy(\un)$. It is established, see~\cite{KLS}: the heights
of the indeterminates $u_1,\dots,u_7$ are constrained to be equal to~$1$, whereas those of $u_8$, $u_9$, $u_{10}$
are not, so there is no concealed parameter here:
\[
\un=(1,1,1,1,1,1,1,N_8,N_9,N_{10}),\quad\text{i.e., }\Par\un=3,
\]
as stated in the Introduction.

The subalgebra of $\fg=\fdy$ generated by $\fg_-$ and $\fg_1$
(the underlined components in table (\ref{1.1})) is, clearly,
isomorphic to $\fbr(3)$. But the new generator of degree 2 and
weight $(0, 0, 2)$ generates, together with $\fbr(3)$, a~larger
algebra.
Below are the dimensions and \textbf{highest} weights of the
components of $\fdy{}^{(1)}$ as $\fdy_0=\fgl(3)$-modules. The terms
of $\fdy$ not contained in $\fdy{}^{(1)}$ are marked in parentheses
in \lq\lq$\dim$\rq\rq column. Recall $(\ref{defdy})$ and
$(\ref{grdy})$; below $i, j, k, l=1, 2, 3$, and $d\omega=0$ is
understood. \scriptsize
\begin{equation}
\label{1.1}
\renewcommand{\arraystretch}{1.4}
\begin{tabular}{|l|l|l|l|}
\hline $\deg$&$\dim$&highest weights of components&basis (for all
possible indices) \cr \hline
$-4$&$3$&$\underline{(-1,-1,-2)}$&$\del_i$\cr
$-3$&$1$&$\underline{(-1,-1,-1)}$&$v^{-1}$\cr
$-2$&$3$&$\underline{(0,-1,-1)}$&$du_iv$\cr
$-1$&$3$&$\underline{(0,0,-1)}$&$du_idu_j$\cr
$0$&$9$&$\underline{(1,0,-1)},\underline{(0,0,0)}$&$u_i\del_j$\cr
$1$&$3$&$\underline{(1,0,0)}$&$u_iv^{-1}$\cr
$2$&$9$&$(2,0,0),\underline{(1,1,0)}$&$u_idu_jv$\cr
$3$&$8$&$(2,1,0)\supset\underline{(1,1,1)}$&$\omega= \sum
f_{ij}du_idu_j\mid \deg f_{ij}=1$\cr $4$&$18$&$(3,1,0),
\underline{(2,1,1)}$&$u_iu_k\del_j$\cr
$5$&$6$&$(3,1,1)$&$u_iu_jv^{-1}$\cr $6$&$18$&$(4,1,1)$,
$(3,2,1)$&$u_iu_jdu_kv$\cr $7$&$15$&$(4,2,1)$&$\omega= \sum
f_{ij}du_idu_j\mid \deg f_{ij}=2$\cr $8$&$21$&$(4,3,1)$, $(4, 2,
2)$&$u_iu_ku_l\del_j$\cr
$9$&$7$&$(4,3,2)$ 
&$u_iu_ku_lv^{-1}$\cr $10$&$21$&$(5,3,2)$,
$(4,4,2)$&$u_iu_ku_ldu_jv$\cr $11$&$15$&$(5,4,2)$&$\omega= \sum
f_{ij}du_idu_j\mid \deg f_{ij}=3$\cr $12$&$18$&$(5,5,2),
(5,4,3)$&$f(u)\del_j\mid \deg f=4$\cr
$13$&$6$&$(5,5,3)$&$f(u)v^{-1}\mid \deg f=4$\cr $14$&$18$&$(6,5,3)$,
$(5,5,4)$&$f(u)du_jv\mid \deg f=4$\cr $15$&$8
(+3)$&$(6,5,4)$&$\omega= \sum f_{ij}du_idu_j\mid \deg f_{ij}=4$\cr
$16$&$9$&$(6,6,4)$, $(6,5,5)$&$f(u)\del_j\mid \deg f=5$\cr
$17$&$3$&$(6,6,5)$&$f(u)v^{-1}\mid \deg f=5$\cr $18$&$9$&$(7,6,5)$,
$(6,6,6)$&$f(u)du_jv\mid \deg f=5$\cr
$19$&$3$&$(7,6,6)$&$\omega= \sum f_{ij}du_idu_j\mid \deg
f_{ij}=5$\cr $20$&$3$&$(7,7,6)$&$f(u)\del_j\mid \deg f=6$\cr
$21$&$1$&$(7,7,7)$&$f(u)v^{-1}\mid \deg f=6$\cr
$22$&$3$&$(8,7,7)$&$f(u)du_jv\mid \deg f=6$\cr \hline
\end{tabular}
\end{equation}
\normalsize

The modules with such highest weights are irreducible if
$\Char\Kee=0$; but since $\Char\Kee=3$, some of these components
are reducible.


\ssec{The big Skryabin algebra $\fg=\fby$}\label{fby} Here $\fg_0=\fgl(3)$
(the weights are given with respect to the $h_i$).
Let $\fg_-$ be realized as follows:
\scriptsize
\begin{equation}\label{realby}
\renewcommand{\arraystretch}{1.4}
\begin{tabular}{|l|l|}
\hline $\fg_{0}\simeq\fgl(V)$&$(-1,0,1)$\quad $u_3 \del_1 + u_3^{(2)}
\del_5 + \left(u_2 u_3-u_4\right) \del_6 -u_2
u_3^{(2)} \del_7$\\
&$(-1,1,0)$\quad $-u_2 \del_1 + u_4 \del_5$\\
&$(0,-1,1)$\quad $-u_3 \del_2 + u_3^{(2)} \del_4 + \left(u_1 u_3 +u_5
\right)\del_6 -u_1
u_3^{(2)} \del_7$;\\
&$h_1=u_1 \del_1 + u_5 \del_5 + u_6 \del_6 + u_7 \del_7,$\\
&$h_2=u_2 \del_2 + u_4 \del_4 + u_6 \del_6 + u_7 \del_7,$\\
&$h_3=u_3 \del_3 + u_4 \del_4 + u_5 \del_5 + u_7 \del_7$;\\
&$(0,1,-1)$\quad $-u_2 \del_3 + u_2^{(2)} \del_4 + \left(u_1 u_2 + u_6
\right)\del_5 -u_1
u_2^{(2)} \del_7$,\\
&$(1,-1,0)$\quad $-u_1 \del_2 + u_5 \del_4$,\\
&$(1,0,-1)$\quad $-u_1 \del_3 + \left(-u_1 u_2 + u_6
\right)\del_4-u_1^{(2)} \del_5 +
u_1^{(2)} u_2 \del_7$\\
\hline $\fg_{-1}\simeq V$&$- \del_1 - u_2\del_6 - u_3 \del_5 + u_4
\del_7,\quad
 - \del_2 + u_3 \del_4 + u_1 \del_6+u_5 \del_7,\ \
- \del_3 + u_6 \del_7$\\
\hline $\fg_{-2}\simeq E^2(V)$&$\del_4 + u_1 \del_7,\quad
  \del_5 + u_2 \del_7,\quad
  \del_6 + u_3 \del_7$\\
\hline
$\fg_{-3}\simeq E^3(V)$&$\del_7$\\
\hline
\end{tabular}
\end{equation}
\normalsize Consider the Cartan-Tanaka-Shchepochkina prolong
$\fby(7;\un):=(\fg_-, \fg_0)_{*,\un}$. As $\fg_0$-module, $\fg_1$ is
a direct sum of two submodules, $\fg'_1$ and $\fg''_1$ with lowest
weights $(0,0, 1)$ and $(-1, 1,1)$, respectively. As algebra,
$\fg_1$ generates 224-dimensional algebra $\fg_+$ of height 14 and
relations up to degree 6 (for comparison: the defining relations of
$\fg_+$ for simple vectorial Lie algebras are of degree 2 (and 3 for
the Hamiltonian series), cf. \cite{GLP}). Observe that even so ugly
and seemingly impossible to use relations are sometimes useful since
they are explicit. 

Let $\fby^{(1)}$ be the algebra generated by $\fg_-$ and $\fg_1$.
Its dimension is 240. The Cartan-Tanaka-Shchepochkina prolong
$(\fg_-, \fg_0)_{*,\un}$ has, however, 4 elements more than
$\fby^{(1)}$: one, of degree 9 and weight
$\underline{\underline{(3,3,3)}}$ and three more elements of degree
12 whose weights are $\underline{\underline{(6,3,3)}}$,
$\underline{\underline{(3,6,3)}}$, and
$\underline{\underline{(3,3,6)}}$. These four elements are outer
derivatives of $(\fg_-, \fg_0)_{*,\un}$; so there are four linearly
independent traces on $(\fg_-, \fg_0)_{*,\un}$ and
$\fby^{(1)}=(\fg_-, \fg_0)_{*,\un}^{(1)}$.

We have
\begin{equation}\label{.27}
[\fg'_1, \fg_{-1}]=\fg_0,\quad [\fg'_1, \fg'_1]=0.
\end{equation}
Let $\fby'$ be the algebra generated by $\fg_-$ and $\fg'_1$. Its
dimension is 19. It is not simple: the part $\fg_{-2}\oplus
\fg_{-3}$ is an ideal. We also have
\begin{equation}\label{.28}
[\fg''_1,\fg_{-1}]=\fsl(3).
\end{equation}
Let $\fby''$ be the algebra generated by $\fg_-$ and $\fg''_1$. It
is the special subalgebra of $\fby$; its dimension is 78. The
element of weight $\underline{\underline{(3,3,3)}}$ is its outer
derivative; together with $\fg_-$, it generates
$\fby'':=\fsby^{(1)}(7;\un)$; the three other outer derivatives of
$\fby$ are also divergence free and belong to the complete prolong
$(\fg_-, \fsl(3))_{*,\un}$.

Here are the dimensions and the \textbf{highest} weights of the
components of $\fby$ (recall that $\fby_0=\fgl(3)$): \scriptsize
\begin{equation}\label{.29}
\renewcommand{\arraystretch}{1.4}
\begin{tabular}{|l|l|l||l|l|l|}
\hline $\deg$&$\dim$&weights of components&$\deg$&$\dim$&weights
of components\cr
\hline
$-3$&$1$&$(-1,-1,-1)$&$6$&$26$&$(4,1,1),(3,2,1),(3,2,1)$\cr
$-2$&$3$&$(0,-1,-1)$&$7$&$24$&$(4,2,1),(3,3,1),(3,2,2)$\cr
$-1$&$3$&$(0,0,-1)$&$8$&$24$&$(4,3,1), (4,2,2),(3,3,2)$\cr
$0$&$9$&$(1,0,-1), (0,0,0)$&$9$&$19$&$(5,2,2),(4,3,2),\underline{\underline{(3,3,3)}}$\cr
$1$&$9$&$(1,1,-1), (1,0,0)$&$10$&$18$&$(5,3,2),(4,3,3)$\cr
$2$&$18$&$(2,1,-1),(1,1,0)$&$11$&$9$&$(5,3,3),(4,4,3)$\cr
$3$&$16$&$(2,1,0), (2,1,0)$&$12$&$11$&$(4,4,4)$, $\underline{\underline{(6,3,3)}}$,
$\underline{\underline{(3,6,3)}}$, $\underline{\underline{(3,3,6)}}$\cr
$4$&$24$&$(3,1,0),(2,2,0),(2,1,1)$&$13$&$3$&$(4,4,4)$\cr
$5$&$24$&$(3,2,0), (3,1,1),(2,2,1)$&$14$&$3$&$(5,5,4)$\cr
 \hline
\end{tabular}
\end{equation}
\normalsize

The dimensions of the respective components of $(\fg_-,
\fsl(3))_{*,\un}$ (same of $\fby''$ only differ in dimensions 6 and
9; the dimensions of outer derivations are paranthesized) are:
\begin{equation}\label{.30}\scriptsize
\renewcommand{\arraystretch}{1.4}
\begin{tabular}{|l|l|l|l|l|l|l|l|l|l|l|l|l|l|}
\hline $\deg$&$-3$&$-2$&$-1$&0&1&2&3&4&5&6&7&8&9\cr \hline
$\dim$&1&3&3&8&6&15&7&15&6&8 (11) &3&3& 0 (1)\cr \hline
\end{tabular}
\end{equation}
The dimension of degree 3 here looks like a~mistake if one compares
with $\deg=3$ line table for $\fby$: but the point is that one
component of weight $(2,1,0)$ contains a~1-dimensional submodule
$(1,1,1)$ while the other component of weight $(2,1,0)$ contains a
submodule of dimension 7.

As $\Zee/2$-graded Lie algebras, these are
\begin{equation}\label{defsby}
\renewcommand{\arraystretch}{1.4}
\begin{array}{l}
\fsby(7;\un):=\fg_\ev \oplus\fg_{\bar 1},\text{~~where
$\fg_\ev\simeq \fsvect(3; \un)$ and $\fg_{\bar 1}\simeq \cO(3;
\un)$}\\
\fsby^{(1)}(7;\un)\simeq \fsvect(3; \un)\oplus \cO'(3; \un),
\end{array}
\end{equation}
where the multiplication of functions is given by
\begin{equation}\label{multsby}
[f, g]=\left(\pderf{f}{u_1}\pderf{g}{u_2}-
\pderf{f}{u_2}\pderf{g}{u_1}\right)\pder{u_3}+cycle(123).
\end{equation}

Hence $\dim \fsby(7;\un)=3^{|\un |+1}+1$ and $\dim\fsby^{(1)}(7;\un
)=3^{|\un |+1}-3$.

Let us now figure out what is $\fby(7;\underline{\widetilde N})$ and
what is its $\underline{\widetilde N}$. It is not difficult to see
that, as space:
\begin{equation}\label{defby}
\renewcommand{\arraystretch}{1.4}
\begin{array}{l}
\fby(\underline{7;\widetilde N}) \simeq \fvect(3; \un)\oplus \cO(3;
\un)_{\Div}\oplus \fvect(3; \un)_{-\Div}\oplus Z^2(3; \un) \
\end{array}
\end{equation}
and we accordingly set:
\begin{equation}\label{grby}
\renewcommand{\arraystretch}{1.4}
\begin{array}{lll}
\deg u^{(\underline{r})}\vvol =2|\underline{r}|-3,&&
\deg u^{(\underline{r})}\del_i =2|\underline{r}|-2,\\
\deg u^{(\underline{r})}\del_iv^{-1} =2|\underline{r}|+1,&&
\deg u^{(\underline{r})}du_idu_j =2|\underline{r}|+4.\\
\end{array}
\end{equation}
The structure of $\fby(7;\underline{\widetilde N})$, as $\fvect(3;
\un)$-module, is remarkable:
\begin{equation}\label{rem}
\renewcommand{\arraystretch}{1.4}
\begin{array}{c}
\fby(7;\underline{\widetilde N}) \simeq
\fvect(3; \un)\oplus Z^2(3; \un)\oplus R,\quad \text{where }\\
0\tto \fvect(3;\un)_{-\Div}\tto R\tto \cO(3;
\un)_{\Div}\tto 0
\end{array}\end{equation}
is an exact sequence.

The multiplication is easy to describe in terms of the vector fields
that constitute the Tanaka-Schepochkina prolongation, but is too
bulky. To describe the multiplication in $\fby(\underline{\widetilde
N})$ in terms of constituents (\ref{rem}), especially the
$\fvect(3;\un)$-action on $R$, observe that the basis of
the complementary module to $\fvect(3;\un)_{-\Div}$ can be
selected canon-ically if we restrict the
$\fvect(3;\un)$-action on $R$ to
$\fsvect(3;\un)$-action or even $\fgl(3)$-action. These
complementary elements will be denoted by \lq\lq$fv$". Let us identify (we usually omit
the wedge product sign between the differentials):
\begin{equation}\label{identif}
du_1du_2:=\del_3v,\quad \del_1\wedge\del_2:=du_3v^2, \text{~~a.s.o.
cyclicly}.
\end{equation}
Then, we have   \tiny
\begin{equation}\label{multby}
\renewcommand{\arraystretch}{1.4}
\begin{tabular}{|l|l|}
\hline
$\fvect\times\fvect\tto\fvect\oplus Z^2$&$Z^2\times\cO_{div}\tto\fvect_{-div}$\\
$[D_1, D_2]_{\fby}= [D_1, D_2]+ d(\Div(D_1)\wedge d(\Div(D_2)$&$[Dv, ``fv"]=-fDv^2$\\

\hline
$\cO_{\Div}\times\cO_{\Div}\tto \fvect$&$\fvect_{-\Div}\times\fvect_{-\Div}\tto Z^2$\\
$[``fv", ``gv"]= (dfdg)v^{2}$&$[D_1v^{-1}, D_2v^{-1}]=[D_1, D_2]v=d(D_1\wedge D_2)v$\\

\hline
$\fvect\times\cO_{div}\tto \fvect_{-\Div}\oplus\cO_{div}$&
$\fvect_{-\Div}\times\cO_{div}\tto\fvect \oplus Z^2$\\
$[D, ``fv"]_{\fby}= \left(D(f)+f\Div D\right)v-df\cdot d(\Div D)$
&$[Dv^{-1}, ``fv"]_{\fby}= fD+ df\cdot d(\Div D)$\\

\hline
$\fvect\times\fvect_{-\Div}\tto \fvect_{-\Div}$&
$\fvect\times Z^2\tto Z^2$\\
Lie derivative
&Lie derivative\\
\hline
$[Z^2,\ Z^2]=0$&$[Z^2,\ \fvect_{-\Div}]=0$.\\
\hline
\end{tabular}
\end{equation}
\normalsize

\sssec{A concealed parameter}\label{concBY} The description
(\ref{defby}) makes an impression that the parameter
$\underline{\widetilde N}$ in $\fby(7;\underline{\widetilde N})$
depends on a~3-dimensional $\un $. Computer experiments show that,
in reality, in the absence of restrictions on $\underline{\widetilde
N}$, the shearing vector depends \textbf{on four parameters, not on
three}:
\begin{equation}\label{N_i}
\underline{\widetilde N}=(1,1,1, N_4, N_5, N_6, N_{7}).
\end{equation}
This is proved in~\cite{KLS}: the heights of $u_1$, $u_2$, $u_3$ are constrained to be equal to~$1$, and the
heights of $u_4,\dots,u_7$ are not. The indeterminates $u_4$, $u_5$, $u_6$ are on equal footing, whereas $u_7$,
the only one of degree~$3$, is not. Looking at the weights of the new elements that appear as we pass
from $N_7=1$ to $N_7=2$, and having observed that
$Z^2\subset\fvect_{\Div}$ in the decomposition \eqref{defby}, we
conjecture that, as a~space,
\begin{equation}\label{BYconj}
\fby(\widetilde N)=\mathop{\oplus}\limits_{0\leq k<
3^{N_7}}\fvect(3;\un)v^k, \text{~~where $\un=(N_4, N_5, N_{6})$},
\text{~~hence $\dim \fby(7;\widetilde N)=3^{|\un|+N_7+1}$}.
\end{equation}
The complete prolong $(\fg_-,\fg_0)_{*,\widetilde N}$ has one more element, as for $\widetilde N=(1,\dots,1)$,
where its dimension is $244=3^5+1$. 
This structure vividly reminds the Lie superalgebra $\fb_\lambda(3)$
of twisted multi-vector fields discovered, it seems, in \cite{Gr0},
see also \cite{LSh}.

\parbegin{Remark} One might think that the existence of the fourth
parameter should be manifest \lq\lq due to symmetry of degree 2 elements of
the irreducible $\fg_0$-module $\fg_{-2}$, see \eqref{realby}". This
argument is wrong as stated: compare with the irreducible
$\fg_0$-module $\fg_{-1}$ in the original construction of the
Melikyan algebra, see \eqref{meme}, where the shearing vector $\un$
with 4 coordinates depends on 2 parameters only. We do not know what
are the conditions the operators of $\fg_0$ should satisfy to ensure
symmetry of the indeterminates from the viewpoint of the critical
dimensions of the shearing vector.
\end{Remark}

\ssec{The middle Skryabin algebra $\fg=\fmy$}\label{fmy}
Here $\fg_0=\fgl(3)$ (the weights are
given with respect to the $H_i$). Let $\fg_-$ be realized by
vector fields as follows:\scriptsize
\begin{equation}\label{realmy}
\renewcommand{\arraystretch}{1.4}
\begin{array}{c}
\renewcommand{\arraystretch}{1.4}
\begin{tabular}{|l|l|l|}
\hline $\fg_{0}$&$(-1,0,1)$\quad $[y_1, y_2]=u_3 \del_1
  + 2 u_3^{(2)}  \del_5+ \left(u_2  u_3
   + 2 u_4\right) \del_6$&$H_1 = u_1 \del_1 + u_5 \del_5
  + u_6 \del_6$\\
&$(-1,1,0)$\quad $y_1=2 u_2 \del_1
    + 2 u_2^{(2)} \del_6
    + u_4 \del_5$&$H_2 = u_2 \del_2 + u_4 \del_4
  + u_6 \del_6$\\
&$(0,-1,1)$\quad $y_2=u_3 \del_2
    + u_3^{(2)} \del_4
    + 2 u_5 \del_6$&$H_3 = u_3 \del_3 + u_4 \del_4
  + u_5 \del_5$ \\
&$(0,1,-1)$\quad $x_2=2 u_2 \del_3
    + 2 u_2^{(2)} \del_4
    + u_6 \del_5$&\\
&$(1,-1,0)$\quad $x_1=2 u_1 \del_2
    + 2 u_1^{(2)} \del_6
    + u_5 \del_4$&\\
&$(1,0,-1)$\quad $[x_1, x_2]=u_1 \del_3+ \left(u_1  u_2 + 2u_6\right) \del_4
    + 2u_1^{(2)} \del_5$&\\
\hline \end{tabular}
\\
\renewcommand{\arraystretch}{1.4}
\begin{tabular}{|l|l|}
\hline $\fg_{-1}$&$\del_1,\; \del_2 + u_1 \del_6,\; \del_3 - u_1
\del_5 + u_2
\del_4$\\
\hline
$\fg_{-2}$&$\del_4,\;\del_5,\;\del_6$\\
\hline \end{tabular}\end{array}
\end{equation}\normalsize
Consider the Cartan-Tanaka-Shchepochkina prolong $(\fg_-, \fg_0)_{*,\un}$. As
$\fg_0$-module, $\fg_1$ is a~direct sum of two submodules,
$\fg'_1$ and $\fg''_1$ with lowest weights $(0,0, 1)$ and $(-1,
1,1)$, respectively. The space $\fg_1$ generates an algebra of
height 11.

Here are the dimensions and \textbf{highest} weights of the
components of $\fmy$, so $\dim\fmy(\One)=162$:\scriptsize
\begin{equation}\label{.31}
\renewcommand{\arraystretch}{1.4}
\begin{tabular}{|l|l|l||l|l|l|}
\hline $\deg$&$\dim$&weights of components&$\deg$&$\dim$&weights
of components\cr \hline
$-2$&$3$&$\underline{(0,-1,-1)}$&$-1$&$3$&$\underline{(0,0,-1)}$\cr
$0$&$9$&$\underline{(1,0,-1)}, \underline{(0,0,0)}$&
$1$&$9$&$(1,1,-1),(1,0,0)$\cr
$2$&$18$&$(2,1,-1),(1,1,0)$&$3$&$18$&$(2,2,-1), (2,1,0)$\cr
$4$&$21$&$(3,1,0), (2,2,0)$&$5$&$21$&$(3,2,0), (3,1,1)$\cr
$6$&$18$&$(4,1,1),(3,2,1)$&
$7$&$18$&$(4,2,1), (3,2,2)$\cr
$8$&$9$&$(4,2,2),(3,3,2)$&
$9$&$9$&$(4,3,2), (3,3,3)$\cr
$10$&$3$&$(4,3,3)$&$11$&$3$&$(4,4,3)$\cr
 \hline
\end{tabular}
\end{equation}
\normalsize

Now, consider partial Cartan-Tanaka-Shchepochkina prolongs $(\fg_-,
\fg_0)_{*,\un}$. As $\fg_0$-module, $\fg_1$ is a~direct sum of two
submodules, $\fg'_1$ and $\fg''_1$ with lowest weights $(0,0, 1)$
and $(-1, 1,1)$, respectively.

Let $\fmy''$ be the algebra generated by $\fg_-$ and $\fg''_1$.
Its negative part is the same as that of $\fmy$,
\begin{equation}\label{.33}
[\fg''_1, \fg_{-1}]=\fsl(3),
\end{equation}
and hence $\fmy''(\un)\simeq\fsmy(\un)$, and $\dim \fsmy(\One)=77$.
Here are the components of $\fsmy$ with weights with respect to the
 matrix units $E_{ii}$ spanning the maximal torus of $\fgl(3)$, for
convenience of comparison with \eqref{.31}: 
\scriptsize
\begin{equation}\label{.51}
\renewcommand{\arraystretch}{1.4}
\begin{tabular}{|l|l|l||l|l|l|}
\hline $\deg$&$\dim$&weights of components&$\deg$&$\dim$&weights of
components\cr \hline
$-2$&$3$&$(0,-1,-1)$&$-1$&$3$&$(0,0,-1)$\cr
$0$&$8$&$(1,0,-1)$& $1$&$6$&$(1,1,-1)$\cr
$2$&$15$&$(2,1,-1)$&$3$&$7$&$(2,1,0)$\cr
$4$&$15$&$(3,1,0)$&$5$&$6$&$(3,1,1)$\cr
$6$&$8$&$(3,2,1)$& $7$&$3$&$(3,2,2)$\cr
$8$&$3$&$(3,3,2)$& &&\cr
 \hline
\end{tabular}
\end{equation}
\normalsize

We have
\begin{equation}
\label{?} [\fg'_1, \fg_{-1}]=\fg_0, \quad [\fg'_1, \fg'_2]=0\ \text{
where }\fg'_2:=[\fg'_1, \fg'_1]. \end{equation}

Let $\fmy'$ be the algebra generated by $\fg_-$ and $\fg'_1$. We
have (given are the \textbf{highest} weights):
\begin{equation}\label{.32}
\renewcommand{\arraystretch}{1.4}
\scriptsize{
\begin{tabular}{|l|l|l||l|l|l|}
\hline $\deg$&$\dim$&weights of components&$\deg$&$\dim$&weights
of components\cr
\hline
$1$&$3$&$(1,0,0)$&$2$&$3$&$(1,1,0)$\cr\hline
\end{tabular}
}
\end{equation}
So $\fmy'\simeq \fo(7)$.

\ssec{Ermolaev algebras $\protect\fg=\protect\fer(3;\tilde\protect\un)$}\label{fer} We have
$\fg_0=\fsl(2)\niplus\fhei(2; 3; \underline{1})$; we realize the
vital components $\fg_i$ for $i=-1, 0, 1$ by vector fields in
indeterminates $x$ as follows (a realization in terms of (\ref{er})
is indicated in parentheses, where $D_i=\del_{u_i}$ to distinguish
from $\del_i=\del_{x_i}$).The weights are given with respect to
$A-B=x_1 \del_1-x_2 \del_2 $ and $B+A=x_1 \del_1 +x_2 \del_2 - 2x_3 \del_3$ from $\fg_0$; in order not to mix the indeterminates
$u_1$ and $u_2$ of realization (\ref{er}), we denote the new three
indeterminates the $x_i$, although they generate the algebra of
divided powers: 
\begin{equation}\label{.34}
\scriptsize
\renewcommand{\arraystretch}{1.4}
\begin{tabular}{|l|l|}
\hline $\fg_{-1}$&$\del_1 ,\quad \del_3 \text{ (this is $\vvol$}),
\quad \del_2$\\
\hline
$\fg_0$&$(-2, 0) : x_2 \del_1 ,$\\
&$(-1, 0) : x_2 \del_3 + x_3 \del_1 \text{ (this is $u_2$})$\\
&$(0, 0) :
A:= x_1 \del_1 - x_3 \del_3 \text{ (this is $u_1D_1$}),$\\
&$ \qquad B:= x_2 \del_2 - x_3 \del_3 \text{ (this is $u_2D_2$}) ,$\\
&$(1, 0) :
   x_1 \del_3 - x_3 \del_2 \text{ (this is $u_1$}),$\\
&$ (2, 0) : x_1 \del_2$\\
\hline
$\fg_1$&$(-3, 1) : - x_2^{(2)} \del_1 ,$\\
&$(-2, 1) :  x_2^{(2)} \del_3 + x_2 x_3 \del_1 ,$\\
&$(-1, 1): [2 \text{ vectors}]:
    x_1 x_2 \del_1 - x_2^{(2)} \del_2,\quad   x_1 x_2 \del_1 - x_2 x_3 \del_3 - x_3^{(2)} \del_1 ,$\\
&$(0, 1) :- x_1 x_2 \del_3 - x_1 x_3 \del_1 +x_2 x_3 \del_2 ,$\\
&$(1, 1): [2 \text{ vectors}]:
   - x_1^{(2)} \del_1 + x_1 x_2 \del_2,\quad  - x_1^{(2)} \del_1 + x_1 x_3 \del_3 - x_3^{(2)} \del_2 ,$\\
&$(2, 1) :
   - x_1^{(2)} \del_3 + x_1 x_3 \del_2 ,$\\
&$(3, 1) :
   - x_1^{(2)} \del_2 $\\
\hline
\end{tabular}
\end{equation}\normalsize
The weights are eigenvalues, i.e., elements of $\Zee/3$: the operator $A+B$ is not homogeneous over~$\Zee$ on, e.g., $x_2\del_3+x_3\del_1$, where it takes the values $3$ and $-3$. So its value is $0$ on all of $\fg_0$ and $1$ on all of $\fg_1$.

The other components of the complete prolong are also computed;
$\fer_1$ is irreducible as a~$\fer_0$-module, it generates the
codimension 1 simple subalgebra $\fer^{(1)}(3;\One)$ of $(\fer_{-1},
\fer_0)_{*,\One}=\mathop{\oplus}\limits_{-1\leq i\leq 4}\fg_i$; the
dimensions of the components of degree 1, 2, 3 are 9, 6 and 2,
respectively; $\dim \fer^{(1)}(\One)=26$.

\sssec{A concealed parameter}\label{concER} We see that $\tilde\un$
for $\fer(3;\tilde\un)$ does not depend on 3 parameters; computer experiments show that
one of indeterminates is critical:
\[
\tilde\un=(\un,1).
\]

\ssec{Frank algebras $\fg=\ffr$} The algebras $\fbr(2; \eps)$,
$\fsp(4)$ and $\ffr(n)$ are partial prolongs with the same
non-positive part as that of $\fk(3; \un)$. The generator of
$\fsp(4)_1$ is $tq$; the generators of $\ffr(n)_1$ are $tp-p^{(2)}q$ and $tq+pq^{(2)}$, 
see eq.~\eqref{.22}, whereas
\begin{equation}\label{.36}
\fbr(2; \eps)_1=\Span(\alpha q^{(2)}p+qt,\; \alpha qp^{(2)}-pt)\;\text{ for
$\alpha=\displaystyle\frac{\eps-1}{\eps+1}$ and $\alpha\neq \pm 1$.}
\end{equation}

The partial prolong of $\left(\mathop{\oplus}\limits_{i\leq 0}
\fk(3; (1, 1, n))_{i}\right)\oplus \ffr(n)_1$ coincides with
$\ffr(n)$ described compo\-nent-wise in \cite{S}.

\section{Defining relations}

The defining relations for the Chevalley generators of $\fbr(2;
\eps)$ are of the same type as Serre relations, but recovered from
the Cartan matrix according to different (as compared with the $p=0$
case), and so far unknown, rules, cf. \cite{GL1}. (Although the
general rules are not known, the answer for $\fbr(2; \eps)$ and
$\fbr(3)$ is now obtained: see relations \eqref{relb2}, and
\eqref{rel1br3}, \eqref{rel2br3}.)

\ssbegin{Theorem}[Presentation of  of $\fbr(2; \eps)$]\label{defrelBr2} For any $\eps\neq 0$, the defining relations between the positive (negative)
generators of $\fbr(2; \eps)$ are
\begin{equation}\label{relb2}
\renewcommand{\arraystretch}{1.4}
\begin{array}{l}
\ad^2_{x_1}(x_2) = 0;\qquad
\ad^3_{x_2}(x_1) = 0.
\end{array}
\end{equation}\end{Theorem}

\ssbegin{Theorem}[Presentation of  of $\fbr(3)$] For $1\fbr(3)_{\pm}$\ , the defining relations are
as follows:\begin{equation}\label{rel1br3} \scriptsize{
\renewcommand{\arraystretch}{1.4}
\begin{array}{l}
[x_1, x_3] = 0;\ \
\ad^2_{x_2}(x_1) = 0,\quad \ad^2_{x_2}(x_3) = 0;\\
\ad^3_{x_3}(x_2) = 0;\ \ {}[[x_3, [x_3, x_2]], [[x_3, [x_2, x_1]],
[x_3, [x_3, x_2]]]] = 0.
\end{array}
}
\end{equation}

For $2\fbr(3)_{\pm}$, we get new gradings of depth $2$ but the
complete prolongs return back the algebra. The relations of degree
$>3$ are different from those in \eqref{rel1br3}:
\scriptsize{\begin{equation}\label{rel2br3}
\renewcommand{\arraystretch}{1.4}
\begin{array}{l}
\ad^3_{x_2}(x_1) = 0,\quad \ad^3_{x_3}(x_2) = 0;\\
{}[[x_1,x_2],[x_3,[x_2,x_3]]]+ [[x_2,x_3,[x_3,[x_1,x_2]]]=0,\quad
{}[[x_2,x_3],[x_2,[x_1,x_2]]]=0.
\end{array}
\end{equation}
}\end{Theorem}

\sssbegin{Remarks} The last non--Serre relation in \eqref{rel1br3} and
the last two relations in \eqref{rel2br3} resemble relations for Lie
superalgebras with Cartan matrix, cf. \cite{GL1}. \end{Remarks}

\ssec{Frank algebras} The relations for $n=1$ are ($x_1=p^{(2)}$): 
\begin{equation}\label{.37} \scriptsize{
\renewcommand{\arraystretch}{1.3}
\begin{array}{ll}
\deg = 1:& [x_1,[x_1,z_1]] = 0,\\
\deg = 2:& [x_1,[x_1,[x_1,z_2]]] = 0,\\
\deg = 3:& [z_1,z_2] = 0, \quad [z_1,[z_1,[x_1,z_1]]] = 0, \quad
[[x_1,z_1],[x_1,[x_1,z_2]]] = 0,\\
\deg = 4:& [z_2,[z_1,[x_1,z_1]]] + [z_2,[x_1,z_2]] = 0,\ \
[[x_1,z_1],[[x_1,z_1],z_2]] + [z_2,[x_1,[x_1,z_2]]] = 0,\\
\deg = 5:& [z_2,[[x_1,z_1],z_2]] = 0,\quad [[x_1,z_2],[[x_1,z_1],z_2]] = 0,\\
\deg = 6:& [z_2,[z_2,[x_1,z_2]]] = 0, \ \
[[x_1,z_2],[z_2,[x_1,z_2]]] = 0,\ \
[[x_1,[x_1,z_2]],[z_2,[x_1,z_2]]] = 0
\end{array}
}
\end{equation}

\subsection*{Acknowledgements} We are thankful to I.~Shchepochkina and S.~Bouarroudj for their help;
DL is thankful to MPIMiS, Leipzig, for financial support and most
creative environment during 2005--07. We thank P.~Zusmanovich and
M.~Kuznetsov for friendly comments and sending us \cite{KD} and
\cite{KuJa}. A half of the updated version was prepared by DL while working at NYUAD; the financial support by the grant AD 065 NYUAD is gratefully acknowledged.

\end{document}